\documentclass[11pt,a4paper,leqno]{article}
\usepackage{a4wide}
\usepackage{latexsym}
\usepackage{amsmath}
\usepackage{amssymb}
\usepackage{stackrel}  
\usepackage{color}
\usepackage{graphicx}
\usepackage[linesnumbered,ruled,vlined]{algorithm2e}

\usepackage{authblk}

\usepackage{hyperref}

\newtheorem{defin}{Definition}
\newtheorem{lemma}{Lemma}
\newtheorem{prop}{Proposition}
\newtheorem{theo}{Theorem}
\newtheorem{corol}{Corollary}
\newtheorem{example}{Example}
\newenvironment{proof}{\medskip\par\noindent{\bf Proof}}{\hfill $\Box$
\medskip\par}

\newcommand{\C}{\mathbb{C}}
\newcommand{\N}{\mathbb{N}}
\newcommand{\R}{\mathbb{R}}
\newcommand{\Z}{\mathbb{Z}}

\begin{document}
\title{On $q-$Gevrey asymptotics for logarithmic type solutions in singularly perturbed $q-$difference-differential equations}

\author[1]{Alberto Lastra}
\author[2]{St\'ephane Malek}
\affil[1]{Universidad de Alcal\'a, Dpto. F\'isica y Matem\'aticas, Alcal\'a de Henares, Madrid, Spain. {\tt alberto.lastra@uah.es}}
\affil[2]{University of Lille, Laboratoire Paul Painlev\'e, Villeneuve d'Ascq cedex, France. {\tt stephane.malek@univ-lille.fr}}


\date{}

\maketitle
\thispagestyle{empty}
{ \small \begin{center}
{\bf Abstract}
\end{center}

A family of singularly perturbed $q-$difference-differential equations under the action of a small complex perturbation parameter is studied. The action of the formal monodromy around the origin is present in the equation, which suggests the construction of holomorphic solutions holding logarithmic terms in both, the formal and the analytic level. We provide both solutions and describe the asymptotic behavior relating them by means of $q-$gevrey asymptotic expansions of some positive order, with respect to the perturbation parameter. 

On the way, the development of a space product of Banach spaces in the Borel plane is needed to provide a fixed point for a coupled system of equations.
\smallskip

}
\bigskip \bigskip

\section{Introduction}

This work is devoted to the study of a family of singularly perturbed $q-$difference-differential equations in the complex domain of the form
\begin{multline}\label{epralintro}
Q(\partial_z)u(t,z,\epsilon)=P(t,z,\epsilon,\sigma_{q;t},\partial_z)u(t,z,\epsilon)+f(t,z,\epsilon)\\
+H\left(\log(\epsilon t),z,\epsilon,u(t,z,\epsilon),\gamma_{\epsilon}^{*}u(t,z,\epsilon)\right).
\end{multline}
Here, the action of a small complex perturbation parameter $\epsilon\in D(0,\epsilon_0)$, for some $\epsilon_0>0$, is studied when considering the unknown function $u=u(t,z,\epsilon)$. 
The element $Q$ stands for a polynomial with complex coefficients. 

The operator $\sigma_{q;t}$ stands for the dilation operator acting on $t$ variable, i.e. $\sigma_{q;t}h(t)=h(qt)$. Here, $q>1$ is a fixed positive number. The previous definition naturally extends to $\sigma_{q;t}^{r}h(t)=h(q^rt)$ for $r\in\mathbb{Q}_{+}$.

The map $P(t,z,\epsilon,\omega_1,\omega_2)$ depends holomorphically on the perturbation parameter $\epsilon$, it is a polynomial in $t$ variable and holomorphic on a horizontal strip in the complex plane with respect to $z$. This function is polynomial with respect to $\omega_1^{1/k}$ and $\omega_2$ variables. Indeed, $k$ will be crucial in the determination of the $q-$Gevrey order relating the formal and the analytic solutions of (\ref{epralintro}). The forcing term $f(t,z,\epsilon)$ is of logarithmic type in the sense that 
\begin{equation}\label{e99}
f(t,z,\epsilon)=f_0(t,z,\epsilon)+f_1(t,z,\epsilon)\frac{\log(\epsilon t)}{\log(q)},
\end{equation}
for some polynomials with respect to $t$ with coefficients being holomorphic functions on some horizontal strip times $D(0,\epsilon_0)$. The concrete form of the forcing term is determined in Assumption (B1) in Section~\ref{secmain}. Finally, $H(v_0,z,\epsilon,v_1,v_2)$ is a linear map with respect to $v_1$, $v_2$, polynomial of degree at most two in $v_0$ and holomorphic on a horizontal strip with respect to $z$ variable and on $D(0,\epsilon_0)$ regarding the perturbation parameter $\epsilon$.

The concrete form of (\ref{epralintro}) is precised in (\ref{epral}), together with the assumptions needed for the main results of the work.

The appearance of the operator $\gamma_{\epsilon}^{*}$ is crucial in the asymptotic behavior of the analytic solutions with respect to $\epsilon$ at the origin. Indeed, the action of the so-called monodromy operator around the origin with respect to $\epsilon$, $\gamma_{\epsilon}^{*}$, described in Section~\ref{secmonoope}, has been proved to modify the form of the analytic solutions of functional equations in the literature, such as~\cite{ma23,ya}. In addition to this, the operator $\gamma_{\epsilon}^{*}$ is used to construct the fundamental solutions to linear systems of differential equations showing an irregular singularity at a point (see Levelt-Turrittin theorem in~\cite{putsinger}, and the reference~\cite{ba2}). The appearance of a monodromy matrix in the previous situation extends the one considered in the regular singular framework which can be obtained by analytically continuing a holomorphic solution defined on a disc out of the origin and turn it counterclockwise around the origin (see~\cite{ilyashyakov}). In a more abstract setting, the monodromy operator in the study of Picard-Vessiot rings has been proved to be important from an algebraic point of view~\cite{putsinger}.

There exists another notion of monodromy used in the study of the meromorphic solutions of $q-$difference equations, introduced in~\cite{rasazh}, which is not considered in the present study.

As a matter of fact, the appearance of the $q-$difference operator $\sigma_{q;t}$ and the monodromy operator $\gamma_{\epsilon}^{*}$ can be reinterpreted as difference-$q-$difference actions on $t$ and $\epsilon$ by the identification of $\gamma_{\epsilon}^{*}$ with the shift mapping of angles $\theta\mapsto \theta+2\pi$ in polar coordinates in the writing $u(t,z,\epsilon)=v(t,z,r,\theta)$, with $\epsilon=re^{i\theta}$ by
$$\gamma_{\epsilon}^{*}u(t,z,\epsilon)=v(t,z,r,\theta+2\pi).$$

In this sense, we refer to~\cite{scsi} and the references therein. In that work, the authors study systems of linear differential and difference equations $Y'(x)=A(x)Y(x)$, together with $\sigma_{q;x}Y(x)=B(x)Y(x)$ or $Y(x+a)=B(x)Y(x)$ for some shift $a$. See also~\cite{romanenko,wang} among other researches in this direction. The study of delay operators and $q-$difference operators is also a point of main interest for applied researchers nowadays. See, for example, the works~\cite{anba,avra,kuku,yaliqi} for the first operators, and~\cite{prrasp,prrasp1} for the second type of operators.

The present work can be seen as a linearized $q-$analog of the problem recently studied by the second author in~\cite{ma23}, where he provides the analytic and formal solutions to a problem of the form
\begin{multline*}
Q(\partial_z)u(t,z,\epsilon)=P(t,z,\epsilon,t\partial_t,\partial_z)u(t,z,\epsilon)+f(t,z,\epsilon)\\
+H\left(\log(\epsilon t),z,\epsilon,\{P_j(\partial_z)u(t,z,\epsilon)\}_{j\in J_1},\{Q_j(\partial_z)\gamma_{\epsilon}^{*}u(t,z,\epsilon)\}_{j\in J_2}\right),
\end{multline*}
(where $\{P_j\}_{j\in J_1}$, $\{Q_j\}_{j\in J_2}$ are finite sets of complex polynomials and $H$ is some nonlinear map in all its arguments) in which the Fuchsian operator $t\partial_t$ is substituted by the $q-$difference operator. The procedure followed in~\cite{ma23} rests on the substitution of the main equation by a decoupled system of equations in the Borel plane, which can be solved recursively as it is a triangularized system of equations. On the contrary, the $q-$analog treated in the present work is more involved, since the coupling of equations in the present study does not allow to come up to a triangular problem. Indeed,  the system is no longer reducible in general. A further advance with respect to the preceding work, is that the general irreducible case is solved, whereas a product Banach space of functions is needed in order to solve the problem. This is the first appearance of such function space in the list of our joint works. We have also included a simplified situation of (\ref{epralintro}) in which a decoupled system arises, and the solution can be obtained by backward substitution of the partial solutions of a triangular system. The last section (Section~\ref{sec10}) is motivated by the more relaxed assumptions needed on the $q-$Gevrey order, being able to solve the main equation under less restrictive geometric conditions. 

We have decided not to state our principal result in the most generality which would have introduced high technical difficulties hiding the main properties of the analytic and formal solutions and their asymptotic relations. In this respect, we have considered functions which remain constant with respect to time variable (see $b_{jk}$ and $c_{\ell}$ in (\ref{epral})), while considering a forcing term of the form (\ref{e99}) instead of the more general 
$$f(t,z,\epsilon)=\sum_{j=0}^{K}f_{j}(t,z,\epsilon)\frac{\log^j(\epsilon t)}{\log^j(q)}.$$ 

The monodromy terms included  in the main equation turn out to have an impact on the $q$-Gevrey order of the asymptotic solution in the general case. In the more restrictive situation in which the associated system is reducible, i.e. triangular, this influence is not so strong, being able to consider less restrictive assumptions. Indeed, the general framework deals with an equation (\ref{epral}) where the $q$-difference operator $\sigma_{q;t}$ appears with rational powers strictly less than 1 in contrast to the usual settings found in the literature where only integer powers arise. In our restrictive reducible case, it is worth noticing that such integer powers are granted. 

The procedure followed in the present study is to search for solutions of (\ref{epralintro}) in the form 
$$u(t,z,\epsilon)=u_0(t,z,\epsilon)+u_1(t,z,\epsilon)\frac{\log(\epsilon t)}{\log(q)},$$
as a natural consequence not only of the form of the forcing term (\ref{e99}) (see Assumption (B1)), but also on the appearance of the monodromy operator within the structure of the equation. This previous motivation yields the pair $\{u_0,u_1\}$ to be a solution of an auxiliary system of equations. At this point, one can search for $u_j$ for $j\in\{0,1\}$ in the form of a so-called $q$-Laplace and inverse Fourier transform, which transforms the previous system into an auxiliary system of equations in two auxiliary unknown functions $\{\omega_0,\omega_1\}$, say
\begin{equation}\label{e129}
\omega_0=\mathfrak{G}_1(\omega_0,\omega_1),\quad \omega_1=\mathfrak{G}_2(\omega_0,\omega_1),
\end{equation}
(see (\ref{e376}), (\ref{e377})).

More precisely, the solution is written in the form
$$u_j(t,z,\epsilon)=\frac{1}{(2\pi)^{1/2}}\frac{k}{\log(q)}\int_{-\infty}^{\infty}\int_{L_d}\frac{\omega_{j}(u,m,\epsilon)}{\Theta_{q^{1/k}}\left(\frac{u}{\epsilon t}\right)}\frac{du}{u}\exp(izm)dm,\quad j\in\{0,1\},$$
for $L_d=[0,\infty)e^{id}$ for some $d\in\R$, and where $\Theta(\cdot)$ stands for Jacobi Theta function (see Section~\ref{secreview}).

The existence of the analytic solution of the auxiliary system (\ref{e129}) needs an accurate description of the geometry and further assumptions on it (see Assumption (D)).  Indeed, strict minorations for the auxiliary polynomial $P_m(\tau)$ (see~(\ref{e434})) need to be provided (see Section~\ref{sec6}). Such tight estimates are needed to compensate the additional terms that come from the presence of logarithmic terms in the problem. Notice that the lower bounds for such a polynomial reached in the works~\cite{drelasmal,ma17} are not sufficient in our new setting.

The solution of the system (\ref{e129}) is attained via a fixed point argument of certain operator acting on a Banach space product. As mentioned above, although the problem enbraces a wide family of functional equations, the appearance of logarithmic terms makes it necessary to adopt several restrictive assumptions on the elements involved in the problem, having a control on its geometry. For that reason, we have decided to include an illustrative concrete example at the end of Section~\ref{sec6}. The first main result of the work (Theorem~\ref{teo1}) describes the form of the analytic solutions to (\ref{epralintro}). The asymptotic behavior of such solutions is also analyzed with respect to the perturbation parameter $\epsilon$ at the origin, through the application of a $q-$analog of the well-known Ramis-Sibuya Theorem (see Theorem ($q-$RS)). Indeed, the $q-$exponential decrease at the intersection of two consecutive sectors in $\epsilon$ uniformly on the rest of the variables (Theorem~\ref{teodif}) is the key point to prove the existence of a formal solution to (\ref{epralintro}) in the form
$$\hat{u}(t,z,\epsilon)= \hat{u}_0(t,z,\epsilon)+\hat{u}_1(t,z,\epsilon)\frac{\log(\epsilon t)}{\log(q)}.$$
In addition to this, for $j\in\{0,1\}$, we prove that $u_j$ admits $\hat{u}_j$ as its $q-$Gevrey asymptotic expansion of some $q-$Gevrey order (see the main result of the present work, Theorem~\ref{teopral}).

We have included a final section, Section~\ref{sec10}, in which the form of the main problem is slightly simplified. Indeed, only one of the coefficients of the monodromy terms in the equation does not appear. This slight simplification leads us to essential simplifications in the geometric assumptions of the problem. More precisely, under these settings, the geometric Assumption (D) is no longer needed. Moreover, the system (\ref{e129}) turns out to be triangular, and therefore much easier to solve. Also, a product Banach space is no longer needed in order to solve the system, but two partial fixed point steps.

The work is structured as follows: In Section~\ref{secpre}, we recall the main facts about the formal monodromy operator, review some integral transforms and their related properties which are involved in the transformation of the main problem into the auxiliary problem. We finally recall the notions of $q-$Gevrey asymptotic expansions and a $q-$Gevrey analog of Ramis-Sibuya Theorem. In Section~\ref{secmain} we state the main problem under study (\ref{epral}) and describe the concrete  assumptions and constructions related to it. After this, we show the strategy we follow to solve the main problem in Section~\ref{secproblemstrategy} and describe the auxiliary Banach spaces of functions involved in the construction of the solution of an auxiliary system in the next section. The analytic solution to the main problem is constructed in Section~\ref{sec7} from the solution of the auxiliary system in Section~\ref{sec6}, where a geometric assumption linked to the solution is described. The formal solution and the asymptotic study relating analytic and formal solution (Theorem~\ref{teopral}) is studied in detail in Section~\ref{sec9}. The work concludes with Section~\ref{sec10}, where the mentioned slight simplification of the main problem is considered.

\vspace{0.3cm}

\textbf{Notation:}

Given $r>0$, $D(0,r)$ denotes the open disc centered at the origin and radius $r$, and $\overline{D}(0,r)$ stands for its closure. Given an open sector $S$ with vertex at the origin, we say that $S_1$ is a subsector of $S$, and denote it by $S_1\prec S$, if $S_1$ is a bounded sector with vertex at the origin and $\overline{S_1}\setminus\{0\}\subseteq S$.

Given a complex Banach space $\mathbb{E}$, we write $\mathbb{E}[[z]]$ for the vector space of formal power series in the variable $z$, with coefficients in $\mathbb{E}$. For every open set $U\subseteq\C$, we write $\mathcal{O}(U,\mathbb{E})$ for the set of holomorphic functions in $U$ with values in $\mathbb{E}$. We adopt the simplified notation $\mathcal{O}(U)$ whenever $\mathbb{E}=\C$. 

Let $q\in\C$. We denote $\sigma_{q;t}$ the dilation operator acting on $t$ variable, i.e. $\sigma_{q;t}f(t)=f(qt)$. This definition is naturally extended to $\sigma_{q;t}^{r}f(t)=f(q^rt)$ for any positive rational number $r$.

\section{Preliminary results}\label{secpre}

In this section, we recall some known results regarding the formal monodromy operator. We also recall the definition and some of the main properties about $q-$Laplace and inverse Fourier transform. The section concludes with a brief review on $q-$Gevrey asymptotic expansions and related results.

\subsection{Formal monodromy operator}\label{secmonoope}
In this subsection, we describe the action of the formal monodromy operator $\gamma_{\epsilon}^{*}$ around the origin. The notion of this operator is stated in Section 3.2~\cite{putsinger}. 

Let us fix $\beta'>0$ and consider the horizontal strip $H_{\beta'}=\{z\in\C:|\hbox{Im}(z)|<\beta'\}$, together with a bounded sector $\mathcal{T}$ with vertex at the origin. We denote $\mathcal{O}_b(\mathcal{T}\times H_{\beta'})$ the Banach space of bounded holomorphic functions on $H_{\beta'}\times\mathcal{T}$ endowed with the sup norm. The next construction can also be generalized in a natural way when substituting $\mathcal{O}_b(\mathcal{T}\times H_{\beta'})$ by any other complex Banach space in $t$ (or $(t,z)$) variable(s). We also set $q>1$.

\begin{defin}
Let $\hat{\mathcal{W}}_1$ be the set of formal power series of the form
\begin{equation}\label{e108}
\hat{u}(t,z,\epsilon)=\hat{u}_0(t,z,\epsilon)+\hat{u}_1(t,z,\epsilon)\frac{\log(\epsilon t)}{\log(q)}
\end{equation}
with $\hat{u}_j\in\mathcal{O}_b(\mathcal{T}\times H_{\beta'})[[\epsilon]]$ for $j=0,1$, and where $\log(\cdot)$ stands for the principal value of the logarithm function defined in $\C\setminus(-\infty,0]$. The formal monodromy operator around the origin in $\C$ with respect to $\epsilon$ is the operator $\gamma_{\epsilon}^{*}$ defined on $\hat{\mathcal{W}}_1$ by
\begin{equation}\label{e125}
\gamma_{\epsilon}^{*}\hat{u}(t,z,\epsilon)=\hat{u}_0(t,z,\epsilon)+\hat{u}_1(t,z,\epsilon)\frac{\log(\epsilon t)+2\pi i}{\log(q)}.
\end{equation}
\end{defin}

\begin{lemma}\label{lema115}
In the situation of the previous definition, it holds that $\gamma_{\epsilon}^{*}(\hat{\mathcal{W}}_1)\subseteq\hat{\mathcal{W}}_1$.
\end{lemma}
\begin{proof}
A rearrangement of the terms in $\gamma_{\epsilon}^{*}\hat{u}(t,z,\epsilon)$ we write
$$\gamma_{\epsilon}^{*}\hat{u}(t,z,\epsilon)=\hat{v}_0(t,z,\epsilon)+\hat{v}_1(t,z,\epsilon)\frac{\log(\epsilon t)}{\log(q)},$$
with $\hat{v}_0(t,z,\epsilon)=\hat{u}_0(t,z,\epsilon)+\frac{2\pi i}{\log(q)}\hat{u}_1(t,z,\epsilon)$ and $\hat{v}_1(t,z,\epsilon)=\hat{u}_1(t,z,\epsilon)$. Observe that $\hat{v}_j\in\mathcal{O}_b(\mathcal{T}\times H_{\beta'})[[\epsilon]]$, for $j=0,1$.
\end{proof}

A natural definition of sum of two elements of $\hat{\mathcal{W}}_1$ and the product with complex numbers provides $\hat{\mathcal{W}}_1$ with the structure of a vector space. It is straight to check the following result. 

\begin{lemma}\label{lema116}
The operator $\gamma_{\epsilon}^{*}:\hat{\mathcal{W}}_1\to\hat{\mathcal{W}}_1$ is a linear endomorphism. 
\end{lemma}

\textbf{Remark:} A more general definition of $\gamma_{\epsilon}^{*}$ can be extended to the set of formal power series of the form
\begin{equation}\label{e126}
\hat{u}(t,z,\epsilon)=\sum_{h=0}^{K}\hat{u}_{h}(t,z,\epsilon)\left(\frac{\log(\epsilon t)}{\log(q)}\right)^{h},
\end{equation}
say $\hat{\mathcal{W}}_K$, with $\hat{u}_h\in\mathcal{O}_b(\mathcal{T}\times H_{\beta'})[[\epsilon]]$. It is straight to check that for every $\hat{u}(t,z,\epsilon)\in\hat{\mathcal{W}}_K$ of the form (\ref{e126}), the formal operator 
$$\gamma_{\epsilon}^{*}\hat{u}(t,z,\epsilon)=\sum_{h=0}^{K}\hat{u}_{h}(t,z,\epsilon)\left(\frac{\log(\epsilon t)}{\log(q)}+\frac{2\pi i}{\log(q)}\right)^{h}$$ 
is such that $\gamma_{\epsilon}^{*}(\hat{\mathcal{W}}_{K})\subseteq \hat{\mathcal{W}}_K$ after a rearrangement of the terms in the previous formula.

The previous definitions in the formal setting can also be adapted to the analytic case, following Section 16 in~\cite{ilyashyakov}. For this purpose, we fix a bounded open sector $\mathcal{E}$ with vertex at the origin, and consider the Banach space of holomorphic and bounded functions defined on $\mathcal{T}\times H_{\beta'}\times \mathcal{E}$, i.e. $\mathcal{O}_b(\mathcal{T}\times H_{\beta'}\times \mathcal{E})$, endowed with the sup norm.

\begin{defin}
Let $\mathcal{W}_1$ be the set of holomorphic functions $u(t,z,\epsilon)$ on the open set $U:=\mathcal{T}\times H_{\beta'}\times\mathcal{E}\setminus \{(t,z,\epsilon)\in\C^{3}: \epsilon t\in(-\infty,0]\}$ of the form
\begin{equation}\label{e108b}
u(t,z,\epsilon)=u_0(t,z,\epsilon)+u_1(t,z,\epsilon)\frac{\log(\epsilon t)}{\log(q)}
\end{equation}
with $u_j\in\mathcal{O}_b(U)$ for $j=0,1$. The action of the formal monodromy operator around the origin in $\C$ with respect to $\epsilon$, $\gamma_{\epsilon}^{*}$, defined on $\mathcal{W}_1$ by
$$\gamma_{\epsilon}^{*}u(t,z,\epsilon)=u_0(t,z,\epsilon)+u_1(t,z,\epsilon)\frac{\log(\epsilon t)+2\pi i}{\log(q)}.$$
\end{defin}

Following Lemma~\ref{lema115} and Lemma~\ref{lema116}, it is straight to check the next result.

\begin{corol}
In the situation of the previous definition, it holds that $\gamma_{\epsilon}^{*}(\mathcal{W}_1)\subseteq \mathcal{W}_1$ and the operator $\gamma_{\epsilon}^{*}:\mathcal{W}_1\to\mathcal{W}_1$ is a linear endomorphism, when endowing $\mathcal{W}_1$ with the usual operations of sum and product.
\end{corol}

Again, the previous definition can be extended to a higher number of terms, as declared in the Remark after Lemma~\ref{lema116}. It is also worth mentioning that the sets $\mathcal{T}\times H_{\beta'}$ (resp. $\mathcal{T}\times H_{\beta'}\times\mathcal{E}$) in the definition of $\hat{\mathcal{W}_1}$ (resp. $\mathcal{W}_1$) can be modified in accordance with other needs. Observe moreover that $\gamma_{\epsilon}^{*}$ coincides with the analytic continuation of a function along a loop at the origin traveled counterclockwise, whenever $u_0$ and $u_1$ are holomorphic functions on some neighborhood of the origin with respect to $\epsilon$.

A direct inspection of the formal operator $\gamma_{\epsilon}^{*}$ allows us to write the terms $\hat{u}_j$, $j=0,1$, of an element $\hat{u}\in\hat{\mathcal{W}}_1$ in the form (\ref{e108}) in terms of the operator $\gamma_{\epsilon}^{*}$ applied on $\hat{u}$. 

\begin{lemma}\label{lema170}
Let $\hat{u}\in\hat{\mathcal{W}}_1$ of the form (\ref{e108}). Then, it holds that
\begin{align*}
\hat{u}_0(t,z,\epsilon)&=\hat{u}(t,z,\epsilon)-\left(\frac{1}{2\pi i}(\gamma_{\epsilon}^{*}-\hbox{id})\hat{u}(t,z,\epsilon)\right)\log(\epsilon t),\\
\hat{u}_1(t,z,\epsilon)&=\frac{\log(q)}{2\pi i}(\gamma_{\epsilon}^{*}-\hbox{id})\hat{u}(t,z,\epsilon),
\end{align*}
where $\hbox{id}$ stands for the identity operator in $\mathcal{O}_{b}(\mathcal{T}\times H_{\beta'})[[\epsilon]]$.
\end{lemma}

Also, the previous representation is valid in the analytic case.

\begin{lemma}\label{lema171}
Let $u\in\mathcal{W}_1$ of the form (\ref{e108b}). Then, it holds that
\begin{align*}
u_0(t,z,\epsilon)&=u(t,z,\epsilon)-\left(\frac{1}{2\pi i}(\gamma_{\epsilon}^{*}-\hbox{id})u(t,z,\epsilon)\right)\log(\epsilon t),\\
u_1(t,z,\epsilon)&=\frac{\log(q)}{2\pi i}(\gamma_{\epsilon}^{*}-\hbox{id})u(t,z,\epsilon),
\end{align*}
where $\hbox{id}$ stands for the identity operator in $\mathcal{O}_{b}(\mathcal{T}\times H_{\beta'}\times\mathcal{E})$.
\end{lemma}

\subsection{Review on some integral transforms and related properties}\label{secreview}

In this subsection, we recall the integral transforms involved in the construction of the analytic solution of the main problem under study in this work, together with some of their main properties which will allow us to state asymptotic results linking them to the formal solution, as the perturbation parameter approaches zero.

Let $k\ge1$ be an integer, and $q>1$ be a real number. The next definition of $q-$Laplace transform has been successfully applied in previous researchs in order to solve functional equations, such as~\cite{drelasmal,ma17}, among others. The proofs of the next results can be found in these references, so we omit them.

Jacobi Theta function of order $k$ is defined by
$$\Theta_{q^{1/k}}(z)=\sum_{p\in\Z}\frac{1}{q^{\frac{p(p-1)}{2k}}}z^p,\qquad z\in\C^{\star}.$$
It turns out to be an analytic function in $\C^{\star}$, with an essential singularity at the origin. The zeros of Jacobi Theta function are given by the elements of the set $\{-q^{m/k}:m\in\Z\}$ as it can be deduced from Jacobi's triple formula. The previous assertion can be refined as follows. 

\begin{lemma}\label{lema194}
Let $\Delta>0$. There exists $C_{q,k}>0$ (which does not depend on $\Delta$) such that
$$|\Theta_{q^{1/k}}(z)|\ge C_{q,k}\Delta\exp\left(\frac{k}{2}\frac{\log^2|z|}{\log(q)}\right)|z|^{1/2},$$
for all $z\in\C^{\star}$ such that $|1+zq^{\frac{m}{k}}|>\Delta$ for every $m\in\Z$. 
\end{lemma} 

\begin{defin}\label{defi194}
Let $\rho>0$ and fix an unbounded sector $S_{d}$ with bisecting direction $d\in\R$ and vertex at the origin, and choose $\gamma\in\R$ such that $L_{\gamma}=\R_{+}e^{i\gamma}\subseteq S_{d}$. Given $f\in\mathcal{O}(D(0,\rho)\cup S_d)$, continuous up to $\overline{D}(0,\rho)$ such that there exist $K,\alpha>0$ and $\delta>1$ with
$$|f(z)|\le K\exp\left(\frac{k}{2}\frac{\log^2(|z|+\delta)}{\log(q)}+\alpha\log(|z|+\delta)\right),$$
for every $z\in S_{d}\cup D(0,\rho)$. We define the $q-$Laplace transform of order $k$ of $f$ along direction $\gamma$ as
$$\mathcal{L}_{q;1/k}^{\gamma}(f(z))(T)=\frac{k}{\log(q)}\int_{L_{\gamma}}\frac{f(u)}{\Theta_{q^{1/k}}\left(\frac{u}{T}\right)}\frac{du}{u}.$$
\end{defin}

The lower estimate of Jacobi Theta function stated in Lemma~\ref{lema194} guarantees the convergence of $q-$Laplace operator acting on functions with appropriate $q-$exponential growth along well-chosen directions.

It is worth remarking that 
\begin{equation}\label{e210}
\mathcal{L}_{q;1/k}^{d}(z^n)(T)=(q^{\frac{1}{k}})^{\frac{n(n-1)}{2}}T^n,
\end{equation}
for every $n\ge0$ and any direction $d\in\R$.

\begin{lemma}\label{lema222}
In the situation of Lemma~\ref{lema194} and given $f$ as in Definition~\ref{defi194}, the $q-$Laplace transform of order $k$ of $f$ defines a holomorphic and bounded function on $\mathcal{R}_{\gamma,\Delta}\cap D(0,r_1)$, for every $0<r_1\le q^{\frac{1/2-\alpha}{k}}/2$, with
$$\mathcal{R}_{\gamma,\Delta}=\left\{T\in\C^{\star}:\left|1+\frac{e^{i\gamma}r}{T}\right|\ge\Delta,\hbox{ for all }r\ge0\right\}.$$
\end{lemma}

\noindent\textbf{Remark:} Direction $\gamma$ may vary among the directions contained in $S_d$ to provide the analytic extension of $q-$Laplace transform of order $k$ of $f$.

\vspace{0.3cm}

An important property of $q-$Laplace transform which will allow us to transform the main equation into a coupled system of auxiliary equations is the following.

\begin{prop}\label{prop218}
Let $f$ be as in Definition~\ref{defi194}, and fix $\tilde{\Delta}>0$. For all $\sigma\ge 0$ and $j\ge0$ it holds that
$$T^{\sigma}\sigma_{q;T}^{j}(\mathcal{L}_{q;1/k}^{\gamma}f(z))(T)=\mathcal{L}_{q;1/k}^{\gamma}\left(\frac{z^{\sigma}}{(q^{1/k})^{\frac{\sigma(\sigma-1)}{2}}}\sigma_{q;z}^{j-\frac{\sigma}{k}}f(z)\right)(T),$$
valid for all $T\in\mathcal{R}_{d,\tilde{\Delta}}\cap D(0,r_1)$, with $0<r_1\le q^{\left(\frac{1}{2}-\alpha\right)/k}/2$.
\end{prop}

We also consider in this work the inverse Fourier transform and related Banach spaces of functions of exponential decrease at infinity, which have also been involved in the solution of functional equations in previous works such as~\cite{lama15,lama}.

\begin{defin}\label{defi218}
Let $\beta>0$ and $\mu>1$. We define the set $E_{(\beta,\mu)}$ of all continuous functions $f:\R\to\C$ such that 
$$\left\|f(m)\right\|_{(\beta,\mu)}=\sup_{m\in\R}(1+|m|)^{\mu}\exp(\beta|m|)|f(m)|<\infty.$$
$E_{(\beta,\mu)}$ is a Banach space endowed with the norm $\left\|\cdot\right\|_{(\beta,\mu)}$.
\end{defin}

\begin{defin}
Let $\beta>0$, $\mu>1$ and fix $f\in E_{(\beta,\mu)}$. The inverse Fourier transform of $f$ is defined by
$$\mathcal{F}^{-1}(f)(x)=\frac{1}{(2\pi)^{1/2}}\int_{-\infty}^{+\infty}f(m)\exp(ixm)dm,\quad x\in\R.$$
\end{defin}

\begin{prop}\label{prop219}
In the previous situation, the following properties hold:
\begin{itemize}
\item[(i)] $\mathcal{F}^{-1}(f)$ can be analytically extended to the set $H_{\beta'}=\{z\in\C:|\hbox{Im}(z)|<\beta'\}$ for every $0<\beta'<\beta$.
\item[(ii)] The function $m\mapsto \phi(m)=imf(m)$ is an element of $E_{(\beta,\mu-1)}$. It holds that $\partial_z\mathcal{F}^{-1}(f)(z)=\mathcal{F}^{-1}(\phi)(z)$.
\item[(iii)] Assume moreover that $g\in E_{(\beta,\mu)}$. The convolution product of $f$ and $g$ given by
$$(f\star g)(m)=\frac{1}{(2\pi)^{1/2}}\int_{-\infty}^{\infty}f(m-m_1)g(m_1)dm_1$$
is an element of $E_{(\beta,\mu)}$. It holds that 
$$\mathcal{F}^{-1}(f)(z)\mathcal{F}^{-1}(g)(z)=\mathcal{F}^{-1}(f\star g)(z),\quad z\in H_{\beta}.$$
\end{itemize}
\end{prop}

\subsection{$q-$Gevrey asymptotic expansions and Ramis-Sibuya theorem}\label{sec23}

In this section we recall the main elements in the theory of holomorphic functions defined on sectors admitting $q-$Gevrey asymptotic expansion, including a $q-$Gevrey version of Ramis-Sibuya theorem.

The definition of $q-$Gevrey asymptotic expansions and related results are slight modifications of those considered in~\cite{lama12}, successfully applied in a different problem in~\cite{ma17}.

In the whole section, $(\mathbb{E},\left\|\cdot\right\|_{\mathbb{E}})$ denotes a complex Banach space. We also assume that $q\in\R$ with $q>1$.

\begin{defin}
Given a bounded sector $\mathcal{E}\subseteq\C^{\star}$ with vertex at the origin, and an integer $k\ge1$, we say that $f\in\mathcal{O}(\mathcal{E},\mathbb{E})$ admits a formal power series $\hat{f}(\epsilon)=\sum_{n\ge0}f_n\epsilon^n\in\mathbb{E}[[\epsilon]]$ as its $q-$Gevrey asymptotic expansion of order $1/k$ if for any given subsector $\tilde{\mathcal{E}}\prec\mathcal{E}$ there exist $C,A>0$ such that 
$$\left\|f(\epsilon)-\sum_{n=0}^{N}f_n\epsilon^n\right\|_{\mathbb{E}}\le CA^{N+1}q^{\frac{N(N+1)}{2k}}|\epsilon|^{N+1},$$
for every $\epsilon\in \tilde{\mathcal{E}}$ and $N\ge0$.
\end{defin}

\noindent\textbf{Remark:} A $q-$Gevrey asymptotic expansion is said to be uniform if the estimates in the previous definition hold for all $\epsilon\in\mathcal{E}$. Observe that any function admitting $q-$Gevrey asymptotic expansion of some order in a sector $\mathcal{E}$ admits uniform $q-$Gevrey asymptotic expansion of the same order in any subsector of $\mathcal{E}$.

\vspace{0.3cm}
 
We recall the classical characterization of asymptotic expansions. Its proof is classical and can be found in Proposition 8, p. 66 of~\cite{ba2}, for example. Hence, we omit its proof.

\begin{prop}\label{prop288}
Let $f\in\mathcal{O}(\mathcal{E},\mathbb{E})$, for some bounded sector $\mathcal{E}$ with vertex at the origin. The following statements are equivalent.
\begin{itemize}
\item[a)] $f$ admits the formal power series $\hat{f}(\epsilon)=\sum_{n\ge0}f_n\epsilon^n$ as its asymptotic expansion in $\mathcal{E}$, i.e. for every subsector $\tilde{\mathcal{E}}\prec\mathcal{E}$ and all $N\ge1$ there exist $C(N,\tilde{\mathcal{E}})>0$ such that 
$$\left\|f(\epsilon)-\sum_{n=0}^{N}f_n\epsilon^n\right\|_{\mathbb{E}}\le C(N,\tilde{\mathcal{E}})|\epsilon|^{N+1},$$
for every $\epsilon\in \tilde{\mathcal{E}}$.
\item[b)] For every integer $n\ge0$ and every $\tilde{\mathcal{E}}\prec\mathcal{E}$ the limit
$$\lim_{\epsilon\to0,\epsilon\in\tilde{\mathcal{E}}}f^{(n)}(\epsilon)$$
exists in $\mathbb{E}$.
\end{itemize}
If one of the previous equivalent statements hold, then the previous limit equals $n!f_n$.
\end{prop}

\textbf{Remark:} In other words, the previous result states that if $f\in\mathcal{O}(\mathcal{E},\mathbb{E})$ admits $\hat{f}(\epsilon)=\sum_{n\ge0}\frac{g_n}{n!}\epsilon^n\in\mathbb{E}[[\epsilon]]$ as its asymptotic expansion in $\mathcal{E}$ (at the origin), then for every $\tilde{\mathcal{E}}\prec\mathcal{E}$ and all $n\ge0$
$$\lim_{\epsilon\to 0,\epsilon\in\tilde{\mathcal{E}}}\left\|f^{(n)}(\epsilon)-g_n\right\|_{\mathbb{E}}=0.$$

\vspace{0.3cm}

The set of functions with null asymptotic expansions coincide with those of $q-$exponential decrease at the origin in the following sense.

\begin{lemma}[Lemma 7,~\cite{ma17}]
Let $\mathcal{E}\subseteq\C^{\star}$ be a bounded sector with vertex at the origin, and $k\ge1$ be an integer number. The following statements are equivalent, for any function $f\in\mathcal{O}(\mathcal{E},\mathbb{E})$:
\begin{itemize}
\item $f$ admits the null formal power series as its $q-$Gevrey asymptotic expansion of order $1/k$ in $\mathcal{E}$.
\item For every $\tilde{\mathcal{E}}\prec\mathcal{E}$, there exist $C>0$ and $K\in\R$ such that 
$$\left\|f(\epsilon)\right\|_{\mathbb{E}}\le C\exp\left(-\frac{k}{2\log(q)}\log^2|\epsilon|\right)|\epsilon|^{K},$$
for $\epsilon\in\tilde{\mathcal{E}}$.
\end{itemize}
\end{lemma}

The cohomological criterion known as Ramis-Sibuya theorem (see~\cite{hssi}, Lemma XI-2-6) has a $q-$analog which can be found in Theorem 25,~\cite{lama13}. See also Theorem (q-RS) in~\cite{ma17} for a detailed proof. Before stating the result, we recall the notion of good covering in $\C^{\star}$.

\begin{defin}\label{defigc}
Let $\zeta\ge2$ be an integer and $(\mathcal{E}_p)_{0\le p\le \zeta-1}$ be a finite family of bounded sectors with vertex at the origin under the following conditions: 
\begin{itemize}
\item[(i)] For every $0\le p\le \zeta-1$ one has that $\mathcal{E}_{p}\cap\mathcal{E}_{p+1}\neq \emptyset$, for all $0\le p\le \zeta-1$ ($\mathcal{E}_{\zeta}:=\mathcal{E}_0$), and the intersection of three of them is empty.
\item[(ii)] $\bigcup_{p=0}^{\zeta-1}\mathcal{E}_{p}$ contains a puntured disc at the origin, and it is contained in $D(0,\epsilon_0)$, for some $\epsilon_0>0$.
\end{itemize}
The family $(\mathcal{E}_{p})_{0\le p\le \zeta-1}$ is know as a good covering in $\C^{\star}$.
\end{defin}

\begin{theo}[q-RS]
Let $\zeta\ge2$ be an integer and $(\mathcal{E}_p)_{0\le p\le \zeta-1}$ be a good covering in $\C^{\star}$. For all $0\le p\le \zeta-1$ we consider $G_p\in\mathcal{O}(\mathcal{E}_p,\mathbb{E})$ and define $\Delta_p=G_{p+1}-G_p\in\mathcal{O}(\mathcal{E}_p\cap\mathcal{E}_{p+1},\mathbb{E})$, where by convention we have put $\mathcal{E}_{\zeta}:=\mathcal{E}_0$ and $G_{\zeta}:=G_0$. Assume that for every $0\le p\le \zeta-1$
\begin{itemize}
\item[a)] $G_p$ is bounded as $\epsilon$ approaches 0 in $\mathcal{E}_p$,
\item[b)] $\Delta_p$ admits uniform null $q-$Gevrey asymptotic expansion of order $1/k$ in $\mathcal{E}_p\cap\mathcal{E}_{p+1}$.
\end{itemize}
Then, there exists $\hat{G}\in\mathbb{E}[[\epsilon]]$ which is the common $q-$Gevrey asymptotic expansion of order $1/k$ of $G_p$ in $\mathcal{E}_p$ for every $0\le p\le \zeta-1$.
\end{theo}

\section{Main problem under study}\label{secmain}

Let $D\ge 2$, $k\ge1$ be integer numbers and $q\in\R$ with $q>1$. 

Let $d_{\ell},\Delta_{\ell}$  for $1\le \ell\le D-1$ and $d_D$ be non negative integers and $\delta_{\ell}$ be a non negative rational number for $1\le \ell\le D-1$. Let $\beta>0$ and fix $0<\beta'<\beta$. We also choose $0<\epsilon_0<1$, and $\mu>1$. These parameters are chosen to satisfy Assumption (A) below. We also assume that $\epsilon_0>0$ is small enough in terms of the elements involved in the problem. The precise value will be precised in this section. We also fix polynomials $R_{\ell}$ for $1\le \ell\le D$ and $Q$ satisfying the second part of Assumption (A) and Assumption (C). 

The forcing term $f(t,z,\epsilon)$ and the coefficients $c_{\ell}(z,\epsilon)$ for $1\le \ell\le D-1$ and $b_{jk}(z,\epsilon)$ for $j,k\in\{0,1\}$ are assumed to be holomorphic functions of certain nature and defined on some domain to be precised (see Assumptions (B1) and (B2)).

The main problem under study in the present work is the equation
\begin{multline}
Q(\partial_z)u(t,z,\epsilon)=(\epsilon t)^{d_D}\sigma_{q;t}^{\frac{d_D}{k}}R_D(\partial_{z})u(t,z,\epsilon)+\sum_{\ell=1}^{D-1}\epsilon^{\Delta_\ell}t^{d_{\ell}}c_{\ell}(z,\epsilon)R_{\ell}(\partial_z)\sigma_{q;t}^{\delta_{\ell}}u(t,z,\epsilon)\\
+f(t,z,\epsilon)+\left(b_{00}(z,\epsilon) +b_{01}(z,\epsilon)\frac{\log(\epsilon t)}{\log(q)}\right)\left[u(t,z,\epsilon)-\left(\frac{1}{2\pi i}(\gamma_{\epsilon}^{*}-\hbox{id})u(t,z,\epsilon)\right)\log(\epsilon t)\right]\\
+\left(b_{10}(z,\epsilon)+b_{11}(z,\epsilon)\frac{\log(\epsilon t)}{\log(q)}\right)\left[\frac{\log(q)}{2\pi i}(\gamma_{\epsilon}^{*}-\hbox{id})u(t,z,\epsilon)\right]\label{epral}.
\end{multline}
We make the following assumptions on the elements involved in the previous equation:

\vspace{0.3cm}

\textbf{Assumption (A):}

$$\Delta_{\ell}>d_{\ell}>k \delta_{\ell}\hbox{ and }d_D\ge k\delta_{\ell},\hbox{ for every }1\le \ell\le D-1.$$

$$\mu>\hbox{deg}(R_{\ell})+1,\qquad 1\le \ell\le D-1.$$

\vspace{0.3cm}

\textbf{Assumption (B1):} The forcing term is of the form
$$f(t,z,\epsilon)=f_0(t,z,\epsilon)+f_1(t,z,\epsilon)\frac{\log(\epsilon t)}{\log(q)},$$
where $f_0,f_1$ are polynomials on their first variable with coefficients being holomorphic functions on $H_{\beta'}\times D(0,\epsilon_0)$. For $h\in\{0,1\}$, we write 
$$f_h(t,z,\epsilon)=F_h(\epsilon t,z,\epsilon)=\sum_{m_h\in\Lambda_h}F_{h,m_h}(z,\epsilon) \left(q^{\frac{1}{k}}\right)^{\frac{m_h(m_h-1)}{2}}(\epsilon t)^{m_h},$$
where $\Lambda_h\subseteq\N_0$ is a finite set. In addition to this, $F_{h,m_h}(z,\epsilon)$ is assumed to be  of the form
$$F_{h,m_h}(z,\epsilon)=\frac{1}{(2\pi)^{1/2}}\int_{-\infty}^{\infty}\tilde{F}_{h,m_h}(m,\epsilon)e^{izm}dm,$$
for some $\R\ni m\mapsto\tilde{F}_{h,m_h}(m,\epsilon)$. The previous function is assumed to be holomorphic for $\epsilon\in D(0,\epsilon_0)$. We write
$$\tilde{F}_{h}(u,m,\epsilon)=\sum_{m_h\in\Lambda_h}\tilde{F}_{h,m_h}(m,\epsilon)u^{m_h}.$$
We assume moreover uniform bounds with respect to the perturbation parameter as follows
$$\sup_{\epsilon\in D(0,\epsilon_0),m\in\R}(1+|m|)^{\mu}e^{\beta|m|}|\tilde{F}_{h,m_h}(m,\epsilon)|\le C_F,$$
for some $C_F>0$, valid for every $h\in\{0,1\}$ and $m_h\in\Lambda_{h}$ (in other words, $m\mapsto\tilde{F}_{h,m_h}\in E_{(\beta,\mu)}$ with uniform holomorphic bounds with respect to the perturbation parameter, see Definition~\ref{defi218}).

This entails that $F_{h,m_h}(z,\epsilon)$ is holomorphic on $H_{\beta'}\times D(0,\epsilon_0)$, for every $0<\beta'<\beta$. Observe from (\ref{e210}) and the previous construction that $F_h$ is built as the $q-$Laplace of order $k$ and inverse Fourier transform of some function, along a direction $d\in\R$. Indeed,
\begin{align*}
F_h(\epsilon t,z,\epsilon)&=\frac{k}{\log(q)}\frac{1}{(2\pi)^{1/2}}\sum_{m_h\in\Lambda_{h}}\int_{L_{d}}\int_{-\infty}^{\infty}\frac{u^{m_h}\tilde{F}_{h,m_h}(m,\epsilon)}{\Theta_{q^{1/k}}\left(\frac{u}{\epsilon t}\right)}e^{izm}dm\frac{du}{u}\\
&=\frac{k}{\log(q)}\frac{1}{(2\pi)^{1/2}}\int_{L_{d}}\int_{-\infty}^{\infty}\frac{\tilde{F}_{h}(u,m,\epsilon)}{\Theta_{q^{1/k}}\left(\frac{u}{\epsilon t}\right)}e^{izm}dm\frac{du}{u}
\end{align*}
for $h\in\{0,1\}$. The expression of $F_h$ actually does not depend on the choice of the direction $d$ since it is a polynomial in $t$.

\vspace{0.3cm}

\textbf{Assumption (B2):} The coefficients $c_\ell(z,\epsilon)$ for $1\le \ell\le D-1$, and $b_{jk}(z,\epsilon)$ for $j,k\in\{0,1\}$ are holomorphic functions defined in $H_{\beta'}\times D(0,\epsilon_0)$, constructed as a inverse Fourier transform. More precisely, there exist $m\mapsto C_{\ell}(m,\epsilon)$ and $m\mapsto \tilde{b}_{jk}(m,\epsilon)$ holomorphic with respect to $\epsilon\in D(0,\epsilon_0)$  such that there exist $\mathcal{C}_{B}>0$ with
\begin{equation}\label{e321}
\sup_{\epsilon\in D(0,\epsilon_0),m\in\R}(1+|m|)^{\mu}e^{\beta|m|}|\tilde{b}_{jk}(m,\epsilon)|\le\mathcal{C}_B,
\end{equation}
for $j,k\in\{0,1\}$ and $\mathcal{C}_{C}>0$ such that
$$\sup_{\epsilon\in D(0,\epsilon_0),m\in\R}(1+|m|)^{\mu}e^{\beta|m|}|C_{\ell}(m,\epsilon)|\le\mathcal{C}_C,$$
for every $1\le \ell\le D-1$ (these functions belong to $E_{(\beta,\mu)}$, see Definition~\ref{defi218} of Section~\ref{secreview}, with uniform bounds with respect to the perturbation parameter). We define
$$c_{\ell}(z,\epsilon)=\frac{1}{(2\pi)^{1/2}}\int_{-\infty}^{\infty}C_{\ell}(m,\epsilon)e^{izm}dm,\quad 1\le \ell\le D-1,$$
$$b_{jk}(z,\epsilon)=\frac{1}{(2\pi)^{1/2}}\int_{-\infty}^{\infty}\tilde{b}_{jk}(m,\epsilon)e^{izm}dm,\quad j,k\in\{0,1\}.$$

Regarding the polynomials $R_{\ell}$ and $Q$, we make the next

\textbf{Assumption (C):}  
$$\hbox{deg}(R_D)=\hbox{deg}(Q).$$
Moreover, there exist $\mathfrak{D}_1,\mathfrak{D}_2>0$ such that 
$$\mathfrak{D}_1\le \inf_{m\in\R}\frac{|Q(im)|}{|R_D(im)|}\le\sup_{m\in\R}\frac{|Q(im)|}{|R_D(im)|}\le \mathfrak{D}_2.$$ 
There exist $\tilde{d}\in\R$ and some small enough $\varsigma>0$, to be determined by the geometric configuration of the problem, which do not depend on $m\in\R$ such that $\arg(Q(im)/R_{D}(im))\in[\tilde{d}-\varsigma,\tilde{d}+\varsigma]$. Observe from the previous condition that $Q(im)\neq 0$ and $R_{D}(im)\neq0$ for $m\in\R$. In addition to this, we assume that 
$$\hbox{deg}(R_{\ell})\le \hbox{deg}(Q)\hbox{ for all }1\le \ell\le D-1.$$

\section{Problem-solving strategy}\label{secproblemstrategy}
In this section, we describe the steps followed to solve the main problem under study, stated in Section~\ref{secmain}. The strategy is divided in different steps modifying the initial equation into an equivalent one, and searching for solutions of a special shape. All the elements involved in the construction of the main problem described in Section~\ref{secmain} are maintained in this section.

First, we search for a formal solution to (\ref{epral}) in the form
\begin{equation}\label{e211}
u(t,z,\epsilon)=U_0(\epsilon t,z,\epsilon)+U_1(\epsilon t,z,\epsilon)\frac{\log(\epsilon t)}{\log(q)},
\end{equation}
motivated by the form of the forcing term in (\ref{epral}) and the action of the monodromy operator. Let us write 
\begin{equation}\label{e284}
u_j(t,z,\epsilon):=U_j(\epsilon t, z,\epsilon),\qquad j\in\{0,1\}.
\end{equation}
By plugging $u(t,z,\epsilon)$ into (\ref{epral}), and taking into account Lemma~\ref{lema170}, one can rewrite (\ref{epral}) in the form
\begin{multline*}
Q(\partial_z)u_0(t,z,\epsilon)+Q(\partial_z)u_1(t,z,\epsilon)\frac{\log(\epsilon t)}{\log(q)}=(\epsilon t)^{d_D}R_{D}(\partial_z)\sigma_{q;t}^{\frac{d_D}{k}}\left[u_0(t,z,\epsilon)+u_1(t,z,\epsilon)\frac{\log(\epsilon t)}{\log(q)}\right]\\
+\sum_{\ell=1}^{D-1}\epsilon^{\Delta_\ell}t^{d_{\ell}}c_{\ell}(z,\epsilon)R_{\ell}(\partial_z)\sigma_{q;t}^{\delta_{\ell}}\left[u_0(t,z,\epsilon)+u_1(t,z,\epsilon)\frac{\log(\epsilon t)}{\log(q)}\right]+f_0(t,z,\epsilon)+f_1(t,z,\epsilon)\frac{\log(\epsilon t)}{\log(q)}\\
+b_{00}(z,\epsilon)u_0(t,z,\epsilon)+b_{10}(z,\epsilon)u_1(t,z,\epsilon)+\left(b_{01}(z,\epsilon)u_0(t,z,\epsilon)+b_{11}(z,\epsilon)u_1(t,z,\epsilon)\right)\frac{\log(\epsilon t)}{\log(q)}
\end{multline*}

which is equivalent to

\begin{multline*}
Q(\partial_z)u_0(t,z,\epsilon)+Q(\partial_z)u_1(t,z,\epsilon)\frac{\log(\epsilon t)}{\log(q)}=(\epsilon t)^{d_D}R_{D}(\partial_z)\left[u_0(q^{\frac{d_D}{k}}t,z,\epsilon)+u_1(q^{\frac{d_D}{k}}t,z,\epsilon)\frac{\log(\epsilon t)}{\log(q)}\right.\\
\left.+\frac{d_D}{k}u_1(q^{\frac{d_D}{k}}t,z,\epsilon)\right]+\sum_{\ell=1}^{D-1}\epsilon^{\Delta_\ell}t^{d_{\ell}}c_{\ell}(z,\epsilon)R_{\ell}(\partial_z)\left[u_0(q^{\delta_{\ell}}t,z,\epsilon)+u_1(q^{\delta_{\ell}}t,z,\epsilon)\frac{\log(\epsilon t)}{\log(q)}+\delta_{\ell}u_1(q^{\delta_{\ell}}t,z,\epsilon)\right]\\
+f_0(t,z,\epsilon)+f_1(t,z,\epsilon)\frac{\log(\epsilon t)}{\log(q)}+b_{00}(z,\epsilon)u_0(t,z,\epsilon)+b_{10}(z,\epsilon)u_1(t,z,\epsilon)\\
+\left(b_{01}(z,\epsilon)u_0(t,z,\epsilon)+b_{11}(z,\epsilon)u_1(t,z,\epsilon)\right)\frac{\log(\epsilon t)}{\log(q)}.
\end{multline*}

The coupling of two equations emerges from the previous equality by identification of the terms with and without the term $\log(\epsilon t)$. More precisely, the previous equation splits into (\ref{e235}),(\ref{e236}):

\begin{multline}
Q(\partial_z)u_0(t,z,\epsilon)=(\epsilon t)^{d_D}R_{D}(\partial_z)\left[u_0(q^{\frac{d_D}{k}}t,z,\epsilon)+\frac{d_D}{k}u_1(q^{\frac{d_D}{k}}t,z,\epsilon)\right]\\
+\sum_{\ell=1}^{D-1}\epsilon^{\Delta_\ell}t^{d_{\ell}}c_{\ell}(z,\epsilon)R_{\ell}(\partial_z)\left[u_0(q^{\delta_{\ell}}t,z,\epsilon)+\delta_{\ell}u_1(q^{\delta_{\ell}}t,z,\epsilon)\right]\\
+f_0(t,z,\epsilon)+b_{00}(z,\epsilon)u_0(t,z,\epsilon)+b_{10}(z,\epsilon)u_1(t,z,\epsilon),\label{e235}
\end{multline}

\begin{multline}
Q(\partial_z)u_1(t,z,\epsilon)=(\epsilon t)^{d_D}R_{D}(\partial_z)u_1(q^{\frac{d_D}{k}}t,z,\epsilon)+\sum_{\ell=1}^{D-1}\epsilon^{\Delta_\ell}t^{d_{\ell}}c_{\ell}(z,\epsilon)R_{\ell}(\partial_z)u_1(q^{\delta_{\ell}}t,z,\epsilon)\\
+f_1(t,z,\epsilon)+b_{01}(z,\epsilon)u_0(t,z,\epsilon)+b_{11}(z,\epsilon)u_1(t,z,\epsilon).\label{e236}
\end{multline}

We recall the expression (\ref{e284}) to rewrite the system of coupled equations (\ref{e235}),(\ref{e236}) in the form

\begin{multline}
Q(\partial_z)U_0(T,z,\epsilon)=T^{d_D}R_{D}(\partial_z)\left[U_0(q^{\frac{d_D}{k}}T,z,\epsilon)+\frac{d_D}{k}U_1(q^{\frac{d_D}{k}}T,z,\epsilon)\right]\\
+\sum_{\ell=1}^{D-1}\epsilon^{\Delta_\ell-d_{\ell}}T^{d_{\ell}}c_{\ell}(z,\epsilon)R_{\ell}(\partial_z)\left[U_0(q^{\delta_{\ell}}T,z,\epsilon)+\delta_{\ell}U_1(q^{\delta_{\ell}}T,z,\epsilon)\right]\\
+F_0(T,z,\epsilon)+b_{00}(z,\epsilon)U_0(T,z,\epsilon)+b_{10}(z,\epsilon)U_1(T,z,\epsilon),\label{e237}
\end{multline}

\begin{multline}
Q(\partial_z)U_1(T,z,\epsilon)=T^{d_D}R_{D}(\partial_z)U_1(q^{\frac{d_D}{k}}T,z,\epsilon)+\sum_{\ell=1}^{D-1}\epsilon^{\Delta_\ell-d_{\ell}}T^{d_{\ell}}c_{\ell}(z,\epsilon)R_{\ell}(\partial_z)U_1(q^{\delta_{\ell}}T,z,\epsilon)\\
+F_1(T,z,\epsilon)+b_{01}(z,\epsilon)U_0(T,z,\epsilon)+b_{11}(z,\epsilon)U_1(T,z,\epsilon).\label{e238}
\end{multline}

At this point, we search for the solution of the coupled system of equations (\ref{e237}),(\ref{e238}) in the form of an inverse Fourier and $q-$Laplace transform (see Section~\ref{secreview}). We write $U_j(T,z,\epsilon)$ in (\ref{e211}) in the form
$$U_j(T,z,\epsilon)=\frac{1}{(2\pi)^{1/2}}\frac{k}{\log(q)}\int_{-\infty}^{\infty}\int_{L_d}\frac{\omega_{j}(u,m,\epsilon)}{\Theta_{q^{1/k}}\left(\frac{u}{T}\right)}\frac{du}{u}\exp(izm)dm,$$
for $j\in\{0,1\}$, and for some direction $d\in\R$ to be precised.

In view of the properties of $q-$Laplace and inverse Fourier transform described in Proposition~\ref{prop218} and Proposition~\ref{prop219}, one has that finding a solution $\{U_0,U_1\}$ of the system (\ref{e237}),(\ref{e238}) is reduced to finding a solution $\{\omega_0,\omega_1\}$ of the following system of convolution equations, (\ref{e376}),(\ref{e377}):

\begin{multline}
Q(im)\omega_0(\tau,m,\epsilon)=\frac{\tau^{d_D}R_{D}(im)}{(q^{\frac{1}{k}})^{\frac{d_D(d_D-1)}{2}}}\omega_0(\tau,m,\epsilon)+\frac{d_D}{k}R_{D}(im)\frac{\tau^{d_D}}{(q^{\frac{1}{k}})^{\frac{d_D(d_D-1)}{2}}}\omega_1(\tau,m,\epsilon)\\
+\sum_{\ell=1}^{D-1}\epsilon^{\Delta_\ell-d_{\ell}}\frac{1}{(2\pi)^{1/2}}\int_{-\infty}^{+\infty}C_{\ell}(m-m_1,\epsilon)R_{\ell}(im_1)\left[\frac{\tau^{d_\ell}}{(q^{\frac{1}{k}})^{\frac{d_{\ell}(d_{\ell}-1)}{2}}}\sigma_{q;\tau}^{\delta_{\ell}-\frac{d_\ell}{k}}\omega_0( \tau,m_1,\epsilon)\right.\\
\hfill\left.+\delta_{\ell}\frac{\tau^{d_\ell}}{(q^{\frac{1}{k}})^{\frac{d_{\ell}(d_{\ell}-1)}{2}}}\sigma_{q;\tau}^{\delta_{\ell}-\frac{d_\ell}{k}}\omega_1( \tau,m_1,\epsilon)\right]dm_1\\
+\tilde{F}_0(\tau,m,\epsilon)+\frac{1}{(2\pi)^{1/2}}\int_{-\infty}^{+\infty}\tilde{b}_{00}(m-m_1,\epsilon)\omega_0(\tau,m_1,\epsilon)+\tilde{b}_{10}(m-m_1,\epsilon)\omega_1(\tau,m_1,\epsilon)dm_1,\label{e376}
\end{multline}

\begin{multline}
Q(im)\omega_1(\tau,m,\epsilon)=\frac{\tau^{d_D}R_{D}(im)}{(q^{\frac{1}{k}})^{\frac{d_D(d_D-1)}{2}}}\omega_1(\tau,m,\epsilon)\\
+\sum_{\ell=1}^{D-1}\epsilon^{\Delta_\ell-d_{\ell}}\frac{1}{(2\pi)^{1/2}}\int_{-\infty}^{+\infty}C_{\ell}(m-m_1,\epsilon)R_{\ell}(im_1)\frac{\tau^{d_\ell}}{(q^{\frac{1}{k}})^{\frac{d_{\ell}(d_{\ell}-1)}{2}}}\sigma_{q;\tau}^{\delta_{\ell}-\frac{d_\ell}{k}}\omega_1( \tau,m_1,\epsilon)dm_1\\
+\tilde{F}_1(\tau,m,\epsilon)+\frac{1}{(2\pi)^{1/2}}\int_{-\infty}^{+\infty}\tilde{b}_{01}(m-m_1,\epsilon)\omega_0(\tau,m_1,\epsilon)+\tilde{b}_{11}(m-m_1,\epsilon)\omega_1(\tau,m_1,\epsilon)dm_1.\label{e377}
\end{multline}

\section{Auxiliary Banach spaces of functions}\label{sec66}

This section is devoted to search adequate Banach spaces of functions in which the solutions of the coupled system of equations (\ref{e376}),(\ref{e377}) can be found, together with the description of the action of some continuous operators acting on such spaces. These properties will allow us to state the existence of analytic solutions of the main equation by means of a fixed point argument. Such functional spaces have already been successfully applied in previous studies such as~\cite{drelasmal,lama12,lama17}, under slight modifications adapted to each concrete problem.

In the whole section, we assume $S_d$ is an open unbounded sector with vertex at the origin and bisecting direction $d\in\R$. We also fix $\rho,\delta>0$ such that $dist(S_d\cup D(0,\rho),-\delta)\ge1$, which is possible as long as $\delta$ is large enough and $S_{d}$ is not bisected by the negative real axis. We also fix $\beta,\mu,\alpha>0$, $q>1$, and an integer number $k\ge1$.  

We omit the proofs of the results included in this section, which can be found in detail in the references mentioned before.

\begin{defin}
The set $\hbox{Exp}_{(k,\beta,\mu,\alpha,\rho)}^{q}$ consists of all complex valued functions $(\tau,m)\mapsto h(\tau,m)$ defined on $(S_d\cup \overline{D}(0,\rho))\times \R$, holomorphic w.r.t. $\tau$ on $S_{d}\cup D(0,\rho)$ such that
$$\left\|h(\tau,m)\right\|_{(k,\beta,\mu,\alpha,\rho)}:=\sup_{\stackrel{\tau\in S_{d}\cup\overline{D}(0,\rho)}{m\in\R}}(1+|m|)^{\mu}e^{\beta |m|}\exp\left(-\frac{k}{2}\frac{\log^2|\tau+\delta|}{\log(q)}-\alpha\log|\tau+\delta|\right)|h(\tau,m)|<\infty.$$
The pair $(\hbox{Exp}_{(k,\beta,\mu,\alpha,\rho)}^{q},\left\|\cdot\right\|_{(k,\beta,\mu,\alpha,\rho)})$ is a Banach space.
\end{defin}

\begin{lemma}\label{lema425}
Let $(\tau,m)\mapsto a(\tau,m)$ be a bounded continuous function defined on $(S_d\cup\overline{D}(0,\rho))\times\R$, holomorphic w.r.t. the first variable on $S_d\cup D(0,\rho)$. Then,
$$\left\|a(\tau,m)h(\tau,m)\right\|_{(k,\beta,\mu,\alpha,\rho)}\le\left(\sup_{\stackrel{\tau\in S_{d}\cup\overline{D}(0,\rho)}{m\in\R}}|a(\tau,m)|\right)\left\|h(\tau,m)\right\|_{(k,\beta,\mu,\alpha,\rho)},$$
for all $h\in \hbox{Exp}_{(k,\beta,\mu,\alpha,\rho)}^{q}$.
\end{lemma}

\begin{prop}\label{prop419}
Let $\gamma_j\ge0$ for $j=1,2,3$ with $\gamma_1+k\gamma_3\ge \gamma_2$. Assume that $a_{\gamma_1}(\tau)$ is a holomorphic function on $S_d\cup D(0,\rho)$, continuous up to $S_d\cup \overline{D}(0,\rho)$ such that 
$$|a_{\gamma_1}(\tau)|\le\frac{1}{(1+|\tau|)^{\gamma_1}},\qquad \tau\in S_d\cup \overline{D}(0,\rho).$$
Then, there exists $C_1=C_1(k,q,\alpha,\{\gamma_j\}_{j=1}^{3})>0$ such that
$$\left\|a_{\gamma_1}(\tau)\tau^{\gamma_2}\sigma_{q;\tau}^{-\gamma_3}h(\tau,m)\right\|_{(k,\beta,\mu,\alpha,\rho)}\le C_1 \left\|h(\tau,m)\right\|_{(k,\beta,\mu,\alpha,\rho)},$$
for all $h\in \hbox{Exp}_{(k,\beta,\mu,\alpha,\rho)}^{q}$.
\end{prop}

In a parallel way, one can consider the Banach space $\hbox{Exp}_{(k,\beta,\mu,\alpha)}^{q}$ in which the set $S_{d}\cup \overline{D}(0,\rho)$ in the definition of $\hbox{Exp}_{(k,\beta,\mu,\alpha,\rho)}^{q}$ is substituted by $S_d$. Proposition~\ref{prop419} reads as follows in this settings.

\begin{corol}\label{coro419}
Let $\gamma_j\ge0$ for $j=1,2,3$ with $\gamma_1+k\gamma_3\ge \gamma_2$. Assume that $a_{\gamma_1}(\tau)$ is a holomorphic function on $S_d$, continuous up to its boundary, such that 
$$|a_{\gamma_1}(\tau)|\le\frac{1}{(1+|\tau|)^{\gamma_1}},\qquad \tau\in S_d.$$
Then, there exists $C_1=C_1(k,q,\alpha,\{\gamma_j\}_{j=1}^{3})>0$ such that
$$\left\|a_{\gamma_1}(\tau)\tau^{\gamma_2}\sigma_{q;\tau}^{-\gamma_3}h(\tau,m)\right\|_{(k,\beta,\mu,\alpha)}\le C_1 \left\|h(\tau,m)\right\|_{(k,\beta,\mu,\alpha)},$$
for all $h\in \hbox{Exp}_{(k,\beta,\mu,\alpha)}^{q}$.
\end{corol}



We extend the definition of the convolution operator acting on the previous Banach space.
\begin{defin}
Let $b,f,g:\R\to\C$ and $h:i\R\to\C$. We define the convolution product
$$f\star^{b,h}g:=b(m)\int_{-\infty}^{\infty}f(m-m_1)h(im_1)g(m_1)dm_1.$$
\end{defin}

This operator coincides with the convolution product $f\star g$ defined in Proposition~\ref{prop219} for $b\equiv \frac{1}{(2\pi)^{1/2}}$ and $h\equiv 1$. 

\begin{prop}\label{prop458}
Let $m\mapsto b(m)$ be a real continuous function such that $|b(m)|\le\frac{1}{|Q_1(m)|}$ for every $m\in\R$ and some polynomial $Q_1$ with complex coefficients. Assume that
$$\hbox{deg}(Q_1)\ge \hbox{deg}(Q_2),\qquad Q_1(im)\neq 0,\quad m\in\R.$$
Let $Q_2$ be a polynomial with complex coefficients, $f\in E_{(\beta,\mu)}$ and $g\in\hbox{Exp}_{(k,\beta,\mu,\alpha,\rho)}^{q}$. In addition to this, assume that $\mu>\hbox{deg}(Q_2)+1$. Then, the function $f\star^{b,Q_2}g$ belongs to $\hbox{Exp}_{(k,\beta,\mu,\alpha,\rho)}^{q}$. Moreover, there exists a positive constant $C_2=C_2(b,Q_2,\mu)$ such that $\left\|f\star^{b,Q_2}g\right\|_{(k,\beta,\mu,\alpha,\rho)}\le C_2 \left\|f\right\|_{(\beta,\mu)} \left\|g\right\|_{(k,\beta,\mu,\alpha,\rho)}$.
\end{prop}
\begin{proof}
We observe there exist $C_{Q_1},C_{Q_2}>0$ such that $|Q_1(m)|\ge C_{Q_1}(1+|m|)^{\hbox{deg}(Q_1)}$ and $|Q_2(im)|\le C_{Q_2}(1+|m|)^{\hbox{deg}(Q_2)}$ for every $m\in\R$.

From the definition of the Banach spaces $E_{(\beta,\mu)}$ and $\hbox{Exp}_{(k,\beta,\mu,\alpha,\rho)}^{q}$ and the convolution operator, one derives
\begin{multline*}
\left\|f\star^{b,Q_2}g\right\|_{(k,\beta,\mu,\alpha,\rho)}\le\sup_{\stackrel{\tau\in S_{d}\cup\overline{D}(0,\rho)}{m\in\R}}(1+|m|)^{\mu}e^{\beta|m|}\exp\left(-\frac{k}{2}\frac{\log^2|\tau+\delta|}{\log(q)}-\alpha\log|\tau+\delta|\right)\\
\times |b(m)|\int_{-\infty}^{+\infty}|f(m-m_1)||Q_2(im_1)||g(\tau,m_1)|dm_1\\
\le \sup_{\stackrel{\tau\in S_{d}\cup\overline{D}(0,\rho)}{m\in\R}}(1+|m|)^{\mu}e^{\beta|m|}\frac{1}{C_{Q_1}(1+|m|)^{\hbox{deg}(Q_1)}}\int_{-\infty}^{+\infty}(1+|m-m_1|)^{\mu}e^{\beta|m-m_1|}|f(m-m_1)|\\
\times C_{Q_2}(1+|m_1|)^{\hbox{deg}(Q_2)}(1+|m_1|)^{\mu}e^{\beta|m_1|}\exp\left(-\frac{k}{2}\frac{\log^2|\tau+\delta|}{\log(q)}-\alpha\log|\tau+\delta|\right)\\
\hfill\times |g(\tau,m_1)|e^{-\beta(|m-m_1|+|m_1|)}\frac{1}{(1+|m-m_1|)^{\mu}}\frac{1}{(1+|m_1|)^{\mu}}dm_1\\
\le \frac{C_{Q_2}}{C_{Q_1}}\left\|f\right\|_{(\beta,\mu)} \left\|g\right\|_{(k,\beta,\mu,\alpha,\rho)}\sup_{m\in\R}(1+|m|)^{\mu-\hbox{deg}(Q_1)}\int_{-\infty}^{+\infty}\frac{dm_1}{(1+|m-m_1|)^{\mu}(1+|m_1|)^{\mu-\hbox{deg}(Q_2)}}.
\end{multline*}
In this last step, we have taken into account that $|m|\le |m-m_1|+|m_1|$ for all $m,m_1\in\R$. We conclude the result by upper estimating the previous integral by a positive constant in virtue of Lemma 2.2~\cite{cota2}, from the hypotheses made.
\end{proof}

As a direct consequence of the previous result with $b\equiv \frac{1}{(2\pi)^{1/2}}$, $Q_1\equiv Q_2\equiv 1$ one has the next result.

\begin{corol}\label{coro486}
Let $f\in E_{(\beta,\mu)}$ and $g\in\hbox{Exp}_{(k,\beta,\mu,\alpha,\rho)}^{q}$. Assume that $\mu>1$. Then, $f\star g \in\hbox{Exp}_{(k,\beta,\mu,\alpha,\rho)}^{q}$ and there exists $C_2(\mu)>0$ such that $\left\|f\star g\right\|_{(k,\beta,\mu,\alpha,\rho)}\le C_2 \left\|f\right\|_{(\beta,\mu)}\left\|g\right\|_{(k,\beta,\mu,\alpha,\rho)}$.
\end{corol}

\section{Analytic solutions of the auxiliary system}\label{sec6}

In this section we prove the existence of analytic solutions to the auxiliary system of coupled equations (\ref{e376}),(\ref{e377}). All the elements and assumptions involved in the statement of the main problem described in Section~\ref{secmain} are maintained in this section. In particular, we take for granted Assumptions (A), (B1), (B2) and (C).

As a first step, we provide lower bounds for the polynomial
\begin{equation}\label{e434}
P_m(\tau)=Q(im)-\frac{1}{(q^{\frac{1}{k}})^{\frac{d_D(d_D-1)}{2}}}R_{D}(im)\tau^{d_D},
\end{equation}
uniformly for $m\in\R$. Let $\{q_{\ell}(m):\ell=0,\ldots,d_{D}-1\}$ denote the set of roots of $\tau\mapsto P_m(\tau)$ for any $m\in\R$. We have that 
\begin{equation}\label{e583}
q_{\ell}(m)=\left(\frac{|Q(im)|q^{\frac{d_D(d_D-1)}{2k}}}{|R_D(im)|}\right)^{1/d_D}\exp
\left(i(\arg(\frac{Q(im)}{R_{D}(im)})\frac{1}{d_D}+\frac{2\pi\ell}{d_D})\right),
\end{equation}
for $0\le \ell\le d_D-1$ and $m\in\R$.
Consider $\varsigma>0$ in Assumption (C) small enough in order that an infinite sector $S_d$ of bisecting direction $d\in\R$ exists, avoiding all the roots of $P_m$, i.e. $S_{d}\cap(\cup_{m\in\R,0\le \ell\le d_D-1}q_{\ell}(m))=\emptyset$. Let $0<\rho<\min\{1,\frac{1}{2}q^{(d_D-1)/2k}\mathfrak{D}_1^{1/d_D}\}$, where $\mathfrak{D}_1$ is given in Assumption (C). We recall that $\delta>0$ is such that $dist(S_d\cup D(0,\rho),-\delta)\ge1$. 

\begin{lemma}\label{lema477}
The following statements hold, concerning $P_m(\tau)$ defined in (\ref{e434}).
\begin{enumerate}
\item[(i)] Let $\tau\in S_d$. There exists $u=u(\tau,m)>0$ and $\theta\not\in 2\pi\Z$ such that 
$$|P_m(\tau)|=|Q(im)||1-u^{d_D}e^{id_D\theta}|,$$ 
for all $m\in\R$.
\item[(ii)] Let $\tau\in D(0,\rho)$. Then,
$$|P_m(\tau)|\ge |Q(im)|\left(1-\frac{1}{2^{d_D}}\right).$$
\end{enumerate}
\end{lemma}
\begin{proof}
The first statement holds in virtue of Assumption (C) and by writting $\tau=ue^{i\theta}q_{\ell}(m)$ for some well chosen $u=u(\tau)>0$, $\theta$ and some root $q_{\ell}(m)$ of $P_m(\tau)$. It holds that
$$|P_m(\tau)|=\left|Q(im)-\frac{\tau^{d_D}R_{D}(im)}{(q^{\frac{1}{k}})^{\frac{d_D(d_D-1)}{2}}}\right|=|Q(im)||1-u^{d_D}e^{id_D\theta}|.$$
For the proof of the second statement, we write $\tau=ue^{i\theta}q_{\ell}(m)$, for some $0\le u\le 1/2$, $\theta\in\R$ and some root $q_{\ell}(m)$ of $P_m(\tau)$. Indeed, observe from the choice of $\rho$ and Assumption (C) that
\begin{multline*}
0<u<\frac{1}{|q_{\ell}(m)|}\rho\le \frac{1}{|q_{\ell}(m)|}\frac{1}{2}q^{\frac{d_D-1}{2k}}\mathfrak{D}_1^{1/d_D}=\left(\frac{|R_D(im)|}{|Q(im)|q^{\frac{d_D(d_D-1)}{2k}}}\right)^{1/d_D}\frac{1}{2}q^{\frac{d_D-1}{2k}}\mathfrak{D}_1^{1/d_D}\\
\le \left(\frac{1}{\mathfrak{D}_1}\right)^{1/d_D}\frac{1}{2}\mathfrak{D}_1^{1/d_D}=\frac{1}{2}.
\end{multline*}
The conclusion follows from analogous estimates as for the first part and the choice of $\rho$. 
\end{proof}

The following is a direct consequence of the previous result.

\begin{corol}\label{coro495}
\begin{itemize}
\item[(i)] There exists $C_{D}>0$ such that $|P_m(\tau)|\ge |Q(im)|C_{D}$ for every $m\in\R$ and all $\tau\in S_d\cup D(0,\rho)$. 
\item[(ii)] Write $\tau=ue^{i\theta}q_{\ell}(m)$ as in the proof of the previous Lemma. Then, there exists $\mathfrak{D}_3>0$ such that 
$$\frac{|1-u^{d_D}e^{id_D\theta}|}{(1+|\tau|)^{d_D}}\ge \mathfrak{D}_3$$
for every $\tau\in S_d$. 
\end{itemize}
\end{corol}
\begin{proof}
The first statement is a straightforward consequence of Lemma~\ref{lema477}. The second statement follows from Assumption (C). Indeed, 
$$\frac{|1-u^{d_D}e^{id_D\theta}|}{(1+|\tau|)^{d_D}}=\frac{|1-u^{d_D}e^{id_D\theta}|}{(1+u|q_{\ell}(m)|)^{d_D}}\ge \frac{|1-u^{d_D}e^{id_D\theta}|}{(1+u\mathfrak{D}_2^{1/d_D}q^{(d_D-1)/2k})^{d_D}}.$$
Observe that the numerator of the last quotient does not vanish for any $u>0$, so given $R>2$ it is lower bounded by a positive constant, say $\mathfrak{D}_{31}$, for all $u\in [0,R]$. For $u>R$, it holds that
\begin{align*}
\frac{|1-u^{d_D}e^{id_D\theta}|}{(1+u\mathfrak{D}_2^{1/d_D}q^{(d_D-1)/2k})^{d_D}}\ge\frac{u^{d_{D}}}{(1+u\mathfrak{D}_2^{1/d_D}q^{(d_D-1)/2k})^{d_D}}\min_{u\ge R}dist(1/u^{d_D},C(0,1))\\
\ge\frac{1}{2}\min_{u\ge R}\frac{u^{d_D}}{(1+u\mathfrak{D}_2^{1/d_D}q^{(d_D-1)/2k})^{d_D}}=:\mathfrak{D}_{32}.
\end{align*}
Take $\mathfrak{D}_{3}=\min\{\mathfrak{D}_{31},\mathfrak{D}_{32}\}$. In the previous estimates, $C(0,1)$ stands for the circle centered at 0 and radius 1.
\end{proof}

We finally consider the following technical assumption, needed for the proof of Proposition~\ref{prop641}

\vspace{0.3cm}

\textbf{Assumption (D):} 
$$\frac{d_D}{k}<\frac{1}{2}\mathfrak{D}_1\min\{C_D,\mathfrak{D}_3\},$$
where $\mathfrak{D}_1$ is determined in Assumption (A) and $\mathfrak{D}_3$ and $C_D$ are  linked to the geometry of the problem, determined in Corollary~\ref{coro495}.

\vspace{0.3cm}

In the following results, we fix $\epsilon\in D(0,\epsilon_0)\setminus\{0\}$ and $\alpha>0$.

\begin{prop}\label{prop536}
Let $1\le \ell\le D-1$. There exists a constant $C_{3,\ell}>0$ (only depending on the parameters involved in the problem, not depending on $\epsilon$) such that the function
$$\mathcal{H}_{\ell}(\omega)=\frac{1}{P_m(\tau)}\frac{1}{(2\pi)^{1/2}}\frac{1}{(q^{\frac{1}{k}})^{\frac{d_{\ell}(d_{\ell}-1)}{2}}}\int_{-\infty}^{+\infty}C_{\ell}(m-m_1,\epsilon)R_{\ell}(im_1)\tau^{d_{\ell}}\sigma_{q;\tau}^{\delta_{\ell}-\frac{d_\ell}{k}}\omega(\tau,m_1)dm_1$$
belongs to $\hbox{Exp}_{(k,\beta,\mu,\alpha,\rho)}^{q}$, provided that $\omega\in\hbox{Exp}_{(k,\beta,\mu,\alpha,\rho)}^{q}$. It holds that 
$$\left\|\mathcal{H}_{\ell}(\omega)\right\|_{(k,\beta,\mu,\alpha,\rho)}\le C_{3,\ell} \left\|\omega\right\|_{(k,\beta,\mu,\alpha,\rho)},$$
for all $\omega\in\hbox{Exp}_{(k,\beta,\mu,\alpha,\rho)}^{q}$.
\end{prop}
\begin{proof}
In view of the different lower bounds obtained for $P_m(\tau)$ when considering $\tau\in S_d$ or $\tau\in D(0,\rho)$ (see Lemma~\ref{lema477}), we divide the proof into two parts. 

First, let $\tau\in S_{d}$ and $m\in\R$. Following the proof of Lemma~\ref{lema477}, we adopt the writing $\tau=ue^{i\theta} q_{\ell}(m)$. We observe from Lemma~\ref{lema477} (i) and Corollary~\ref{coro495} (ii) that
$$\frac{1}{|P_{m}(\tau)|}=\frac{1}{|Q(im)||1-u^{d_D}e^{id_D\theta}|}\le \frac{1}{|Q(im)|\mathfrak{D}_3(1+|\tau|)^{d_D}},$$
for all $m\in\R$. Observe from the previous assertion that the function $\tau\mapsto Q(im)/P_{m}(\tau)$ admits uniform upper bounds for $m\in\R$. Due to Assumption (A) holds, we may fix $m\in\R$ and apply Corollary~\ref{coro419} to the function
\begin{equation}\label{e542}
(\tau,m_1)\mapsto H(\tau,m,m_1)=\frac{Q(im)}{P_{m}(\tau)}\tau^{d_{\ell}}\sigma_{q;\tau}^{\delta_{\ell}-\frac{d_{\ell}}{k}}\omega(\tau,m_1).
\end{equation}
On the other hand, for $\tau\in \overline{D}(0,\rho)$, one derives from Lemma~\ref{lema477} (ii) that the function (\ref{e542}) satisfies
\begin{multline*}
|H(\tau,m,m_1)|(1+|m_1|)^{\mu}e^{\beta|m_1|}\exp\left(-\frac{k}{2}\frac{\log^2|\tau+\delta|}{\log(q)}-\alpha\log|\tau+\delta|\right)\\
\le \frac{1}{1-\frac{1}{2^{d_D}}}\rho^{d_{\ell}}\exp\left(\frac{k}{2}\frac{\log^2|q^{\delta_{\ell}-d_{\ell}/k}\tau+\delta|-\log^2|\tau+\delta|}{\log(q)}+\alpha(\log|q^{\delta_{\ell}-d_{\ell}/k}\tau+\delta|-\log|\tau+\delta|)\right)\\
\times \left[(1+|m_1|)^{\mu}e^{\beta|m_1|}\exp\left(-\frac{k}{2}\frac{\log^2|q^{\delta_{\ell}-d_{\ell}/k}\tau+\delta|}{\log(q)}-\alpha\log|q^{\delta_{\ell}-d_{\ell}/k}\tau+\delta|\right)\omega(q^{\delta_{\ell}-d_{\ell}/k}\tau,m_1)\right]\\
\le\frac{1}{1-\frac{1}{2^{d_D}}}\rho^{d_{\ell}}\exp\left(\mathfrak{D}_{4}\right)\left\|\omega\right\|_{(k,\beta,\mu,\alpha,\rho)}
\end{multline*}
for some $\mathfrak{D}_4>0$, valid for all $m,m_1\in\R$. Here,
\begin{multline*}
\max_{\tau\in \overline{D}(0,\rho)} \frac{k}{2}\frac{\log^2|q^{\delta_{\ell}-d_{\ell}/k}\tau+\delta|-\log^2|\tau+\delta|}{\log(q)}+\alpha(\log|q^{\delta_{\ell}-d_{\ell}/k}\tau+\delta|-\log|\tau+\delta|)\\
\le \frac{k}{2\log(q)}(\log^2(\rho+\delta)+\alpha\log(\rho+\delta))=:\mathfrak{D}_4
\end{multline*}
The proof is concluded by applying Proposition~\ref{prop458} with $Q_1:=Q$, $Q_2:=R_{\ell}$, $f:=C_{\ell}$ and $g:=H$, which can be applied regarding the restrictions appearing in Assumption (A). We finally observe that one can choose 
$$C_{3,\ell}:=\max\left\{\frac{1}{1-\frac{1}{2^{d_D}}}\rho^{d_\ell}\exp(\mathfrak{D}_{4}),C_1\right\}C_2\mathcal{C}_{C},$$
where $\mathcal{C}_{C}$ is the constant appearing in Assumption (B2), $C_1$ is given in Corollary~\ref{coro419} and $C_2$ is determined in Proposition~\ref{prop458}.
\end{proof}

\begin{prop}\label{prop565} The function 
$$\mathcal{H}_{P}(\omega_1)=\frac{1}{(q^{\frac{1}{k}})^{\frac{d_D(d_D-1)}{2}}}\frac{R_{D}(im)}{P_m(\tau)}\frac{d_D}{k}\tau^{d_D}\omega_1(\tau,m)$$
belongs to $\hbox{Exp}^{q}_{(k,\beta,\mu,\alpha,\rho)}$, provided that $\omega_1\in\hbox{Exp}^{q}_{(k,\beta,\mu,\alpha,\rho)}$. Moreover, one has that 
$$\left\|\mathcal{H}_{P}(\omega_1)\right\|_{(k,\beta,\mu,\alpha,\rho)}\le\frac{d_D}{k}\frac{1}{\mathfrak{D}_1}\max\left\{\frac{\rho^{d_D}}{C_D},\mathfrak{D}_3^{-1}\right\}\left\|\omega_1\right\|_{(k,\beta,\mu,\alpha,\rho)}\le \frac{d_D}{k}\frac{1}{\mathfrak{D}_1}\max\{\frac{1}{C_D},\frac{1}{\mathfrak{D}_3}\}\left\|\omega_1\right\|_{(k,\beta,\mu,\alpha,\rho)} .$$
\end{prop}
\begin{proof}
We divide the proof into two parts. We first choose $\tau\in \overline{D}(0,\rho)$ and $m\in\R$. Then, in view of Corollary~\ref{coro495} (i), and Assumption (C), it holds that
\begin{multline*}
\left|\frac{R_{D}(im)}{P_m(\tau)}\frac{d_D}{k}\tau^{d_D}\omega_1(\tau,m)\right|(1+|m|)^{\mu}e^{\beta|m|}\exp\left(-\frac{k}{2}\frac{\log^2|\tau+\delta|}{\log(q)}-\alpha\log|\tau+\delta|\right)\\
\le \frac{|R_{D}(im)|}{|Q(im)|C_D}\frac{d_D}{k}\rho^{d_D}|\omega_1(\tau,m)|(1+|m|)^{\mu}e^{\beta|m|}\exp\left(-\frac{k}{2}\frac{\log^2|\tau+\delta|}{\log(q)}-\alpha\log|\tau+\delta|\right)\\
\le \frac{d_D}{k}\frac{1}{\mathfrak{D}_1}\frac{\rho^{d_D}}{C_D}\left\|\omega_1\right\|_{(k,\beta,\mu,\alpha,\rho)}.
\end{multline*}
On the other hand, for all $\tau\in S_d$ and $m\in\R$, regarding Corollary~\ref{coro495} (ii) , Lemma~\ref{lema477}, and Assumption (C), one has that
\begin{multline*}
\left|\frac{R_{D}(im)}{P_m(\tau)}\frac{d_D}{k}\tau^{d_D}\omega_1(\tau,m)\right|(1+|m|)^{\mu}e^{\beta|m|}\exp\left(-\frac{k}{2}\frac{\log^2|\tau+\delta|}{\log(q)}-\alpha\log|\tau+\delta|\right)\\
\le \frac{|R_{D}(im)|}{|Q(im)|}\frac{1}{\mathfrak{D}_3}\frac{d_D}{k}\frac{|\tau|^{d_{D}}}{(1+|\tau|)^{d_D}}|\omega_1(\tau,m)|(1+|m|)^{\mu}e^{\beta|m|}\exp\left(-\frac{k}{2}\frac{\log^2|\tau+\delta|}{\log(q)}-\alpha\log|\tau+\delta|\right)\\
\le \frac{d_D}{k}\frac{1}{\mathfrak{D}_1\mathfrak{D}_3}\left\|\omega_1\right\|_{(k,\beta,\mu,\alpha,\rho)}.
\end{multline*}

The previous bounds yield
$$\left\|\frac{R_{D}(im)}{P_m(\tau)}\frac{d_D}{k}\tau^{d_D}\omega_1(\tau,m)\right\|_{(k,\beta,\mu,\alpha,\rho)}\le\frac{d_D}{k}\frac{1}{\mathfrak{D}_1}\max\left\{\frac{\rho^{d_D}}{C_D},\frac{1}{\mathfrak{D}_3}\right\}\left\|\omega_1\right\|_{(k,\beta,\mu,\alpha,\rho)},$$
valid for all $\tau\in\overline{D}(0,\rho)\cup S_d$ and $m\in\R$. We conclude the proof by taking into account the choice $0<\rho\le 1$. 
\end{proof}

\begin{prop}\label{prop598}
There exists $\tilde{C}_{F}>0$ (only depending on the parameters involved in the problem, not depending on $\epsilon$) such that for $h\in\{0,1\}$, the function $\tilde{F}_{h}(\tau,m,\epsilon)/P_m(\tau)$ belongs to $\hbox{Exp}^{q}_{(k,\beta,\mu,\alpha,\rho)}$ with
$$\left\|\frac{\tilde{F}_{h}(\tau,m,\epsilon)}{P_m(\tau)}\right\|_{(k,\beta,\mu,\alpha,\rho)}\le \tilde{C}_F.$$
\end{prop}
\begin{proof}
Let $h\in\{0,1\}$. In view of Assumption (C), Lemma~\ref{lema477} and Corollary~\ref{coro495} one arrives at
\begin{multline*}
\left\|\frac{\tilde{F}_h(\tau,m,\epsilon)}{P_m(\tau)}\right\|_{(k,\beta,\mu,\alpha,\rho)}\le\sup_{\tau\in S_d\cup\overline{D}(0,\rho),m\in\R}\sum_{m_h\in\Lambda_h}|\tilde{F}_{h,m_h}(m,\epsilon)|(1+|m|)^{\mu}e^{\beta|m|}\frac{|\tau|^{m_h}}{|P_m(\tau)|}\\
\hfill\times\exp\left(-\frac{k}{2}\frac{\log^2|\tau+\delta|}{\log(q)}-\alpha\log|\tau+\delta|\right)\\
\le \sup_{\tau\in S_d\cup\overline{D}(0,\rho),m\in\R}\sum_{m_h\in\Lambda_h}\left\|\tilde{F}_{h,m_h}\right\|_{(\beta,\mu)}\frac{|\tau|^{m_h}}{|P_m(\tau)|}\exp\left(-\frac{k}{2}\frac{\log^2|\tau+\delta|}{\log(q)}-\alpha\log|\tau+\delta|\right)\\
\le C_F\sup_{m\in\R}\left\{\frac{1}{|Q(im)|}\right\}\sum_{m_h\in\Lambda_h}\sup\left\{\frac{1}{1-\frac{1}{2^{d_D}}}\sup_{\tau\in\overline{D}(0,\rho)}|\tau|^{m_h}\exp(-\alpha\log|\tau+\delta|),\right.\\
\hfill\left.\frac{1}{\mathfrak{D}_3}\sup_{\tau\in S_d}\left(\frac{|\tau|^{m_h}}{(1+|\tau|)^{d_D}}\exp\left(-\frac{k}{2}\frac{\log^2|\tau+\delta|}{\log(q)}-\alpha\log|\tau+\delta|\right)\right)\right\}\\
\le C_F\sup_{m\in\R}\left\{\frac{1}{|Q(im)|}\right\}\sum_{m_h\in\Lambda_h}\sup\left\{\frac{1}{1-\frac{1}{2^{d_D}}}\max_{\tau\in D(0,\rho)}|\tau|^{m_h}e^{-\alpha\log|\tau+\delta|},\right.\\
\hfill\left. \frac{1}{\mathfrak{D}_3}\sup_{x>0}\left((x+\delta)^{m_h}\exp(-\frac{k}{2\log(q)}\log^2(x)-\alpha\log(x))\right)\right\}.
\end{multline*}
We observe that the last expression is upper bounded to conclude the result.
\end{proof}

Let us consider the product Banach space $\left(\hbox{Exp}_{(k,\beta,\mu,\alpha,\rho)}^{q}\right)^2$ endowed with the norm 
$$\left\|\begin{pmatrix} f\\g\end{pmatrix}\right\|:=\max\left\{\left\|f\right\|_{(k,\beta,\mu,\alpha,\rho)},\left\|g\right\|_{(k,\beta,\mu,\alpha,\rho)}\right\}.$$
We work with the operator 
$$\mathcal{H}_{\epsilon}: \left(\hbox{Exp}_{(k,\beta,\mu,\alpha,\rho)}^{q}\right)^2\to \left(\hbox{Exp}_{(k,\beta,\mu,\alpha,\rho)}^{q}\right)^2$$
defined by
$$
\mathcal{H}_{\epsilon}\begin{pmatrix} \omega_0(\tau,m)\\\omega_1(\tau,m)\end{pmatrix}=
\begin{pmatrix}\displaystyle
\mathcal{H}_{P}(\omega_1)+\sum_{\ell=1}^{D-1}\epsilon^{\Delta_\ell-d_{\ell}}(\mathcal{H}_{\ell}(\omega_0)+\delta_{\ell}\mathcal{H}_{\ell}(\omega_1)) +\frac{\tilde{F}_0(\tau,m,\epsilon)}{P_m(\tau)}+\frac{1}{P_m(\tau)}\sum_{j=0}^{1}\tilde{b}_{j0}\star\omega_j\\
\displaystyle 
\sum_{\ell=1}^{D-1}\epsilon^{\Delta_\ell-d_{\ell}}\mathcal{H}_{\ell}(\omega_1)+\frac{\tilde{F}_1(\tau,m,\epsilon)}{P_m(\tau)}+\frac{1}{P_m(\tau)}\sum_{j=0}^{1}\tilde{b}_{j1}\star\omega_j
\end{pmatrix},
$$
where $\mathcal{H}_{P}(\cdot)$ is given in Proposition~\ref{prop565}, $\mathcal{H}_{\ell}(\cdot)$ for $1\le \ell\le D-1$ is defined in Proposition~\ref{prop536}, and the convolution product is stated in Proposition~\ref{prop219}.

\begin{prop}\label{prop641}
There exist $\epsilon_0,\varpi,\varsigma_b>0$  such that for every $\epsilon\in D(0,\epsilon_0)\setminus\{0\}$ and $\mathcal{C}_{B}\le \varsigma_b$ (recall that $\mathcal{C}_{B}$ is determined by (\ref{e321})), the operator $\mathcal{H}_{\epsilon}$ satisfies that 
$$\mathcal{H}_{\epsilon}(\overline{B}(0,\varpi)\times\overline{B}(0,\varpi))\subseteq \left(\overline{B}(0,\varpi)\times \overline{B}(0,\varpi)\right),$$
where $\overline{B}(0,\varpi)\subseteq \hbox{Exp}_{(k,\beta,\mu,\alpha,\rho)}^{q}$ stands for the closed disc of radius $\varpi$ centered at the origin in $\hbox{Exp}_{(k,\beta,\mu,\alpha,\rho)}^{q}$. In addition to this, for every $\omega_{jk}\in \overline{B}(0,\varpi) $, for $(j,k)\in\{0,1\}\times\{1,2\}$ it holds that
$$\left\|\mathcal{H}_{\epsilon}\begin{pmatrix} \omega_{01}(\tau,m)\\\omega_{11}(\tau,m)\end{pmatrix}-\mathcal{H}_{\epsilon}\begin{pmatrix} \omega_{02}(\tau,m)\\\omega_{12}(\tau,m)\end{pmatrix}\right\|\le \frac{1}{2}\left\|\begin{pmatrix} \omega_{01}(\tau,m)\\\omega_{11}(\tau,m)\end{pmatrix}-\begin{pmatrix} \omega_{02}(\tau,m)\\\omega_{12}(\tau,m)\end{pmatrix}\right\|.$$
\end{prop}
\begin{proof}
For the first part of the proof, let us consider $\varpi>0$ such that $\tilde{C}_{F}\le \varpi/2$ and choose $\epsilon_0>0$ and $\varsigma_b>0$ such that
\begin{equation}\label{e636}
\frac{d_D}{k}\frac{1}{\mathfrak{D}_1}\max\{\frac{1}{C_D},\frac{1}{\mathfrak{D}_3}\}+\sum_{\ell=1}^{D-1}\epsilon_0^{\Delta_{\ell}-d_\ell}C_{3,\ell}(1+\delta_\ell)+\sup_{m\in\R}\left(\frac{1}{|Q(im)|}\right)\frac{2}{C_D}\varsigma_{b}C_2\le\frac{1}{2},
\end{equation}
where $C_2$ is the constant determined in Corollary~\ref{coro486}. Observe that Assumption (D) is needed in order that (\ref{e636}) holds for adequate $\epsilon_0,\varsigma_b$.

Let $\omega_{j}\in \overline{B}(0,\varpi)$ for $j\in\{0,1\}$. We apply Lemma~\ref{lema425}, Proposition~\ref{prop536}, Proposition~\ref{prop565}, Assumptions (B1) and (B2), together with Corollary~\ref{coro486}, Corollary~\ref{coro495} and Proposition~\ref{prop598} to obtain that
\begin{multline}
\left\|\mathcal{H}_{P}(\omega_1)+\sum_{\ell=1}^{D-1}\epsilon^{\Delta_\ell-d_{\ell}}(\mathcal{H}_{\ell}(\omega_0)+\delta_{\ell}\mathcal{H}_{\ell}(\omega_1)) +\frac{\tilde{F}_0(\tau,m,\epsilon)}{P_m(\tau)}+\frac{1}{P_m(\tau)}\sum_{j=0}^{1}\tilde{b}_{j0}\star\omega_j\right\|_{(k,\beta,\mu,\alpha,\rho)}\\
\le \frac{d_D}{k}\frac{1}{\mathfrak{D}_1}\max\{\frac{1}{C_D},\frac{1}{\mathfrak{D}_3}\}\left\|\omega_1\right\|_{(k,\beta,\mu,\alpha,\rho)}+\sum_{\ell=1}^{D-1}\epsilon_0^{\Delta_{\ell}-d_\ell}C_{3,\ell}\left(\left\|\omega_0\right\|_{(k,\beta,\mu,\alpha,\rho)}+\delta_{\ell}\left\|\omega_1\right\|_{(k,\beta,\mu,\alpha,\rho)}\right)\\
\hfill +\tilde{C}_F+\sup_{m\in\R}\left(\frac{1}{|Q(im)|}\right)\frac{1}{C_D}\mathcal{C}_{B}C_2\left(\left\|\omega_0\right\|_{(k,\beta,\mu,\alpha,\rho)}+\left\|\omega_1\right\|_{(k,\beta,\mu,\alpha,\rho)}\right)\\
\le \left[\frac{d_D}{k}\frac{1}{\mathfrak{D}_1}\max\{\frac{1}{C_D},\frac{1}{\mathfrak{D}_3}\}+\sum_{\ell=1}^{D-1}\epsilon_0^{\Delta_{\ell}-d_\ell}C_{3,\ell}(1+\delta_\ell)+\sup_{m\in\R}\left(\frac{1}{|Q(im)|}\right)\frac{2}{C_D}\mathcal{C}_{B}C_2\right]\varpi+\tilde{C}_F\\
\hfill \le\frac{\varpi}{2}+\frac{\varpi}{2}=\varpi.\label{e647}
\end{multline}
The choice in (\ref{e636}) allows us to conclude. An analogous reasoning yields
\begin{multline}
\left\|\sum_{\ell=1}^{D-1}\epsilon^{\Delta_\ell-d_{\ell}}\mathcal{H}_{\ell}(\omega_1)+\frac{\tilde{F}_1(\tau,m,\epsilon)}{P_m(\tau)}+\frac{1}{P_m(\tau)}\sum_{j=0}^{1}\tilde{b}_{j1}\star\omega_j\right\|_{(k,\beta,\mu,\alpha,\rho)}\\
\le \sum_{\ell=1}^{D-1}\epsilon_0^{\Delta_{\ell}-d_\ell}C_{3,\ell}\left\|\omega_1\right\|_{(k,\beta,\mu,\alpha,\rho)}+ \tilde{C}_F+\sup_{m\in\R}\left(\frac{1}{|Q(im)|}\right)\frac{1}{C_D}\mathcal{C}_{B}C_2\left(\left\|\omega_0\right\|_{(k,\beta,\mu,\alpha,\rho)}+\left\|\omega_1\right\|_{(k,\beta,\mu,\alpha,\rho)}\right)\\
\le \left[\sum_{\ell=1}^{D-1}\epsilon_0^{\Delta_{\ell}-d_\ell}C_{3,\ell}+\sup_{m\in\R}\left(\frac{1}{|Q(im)|}\right)\frac{2}{C_D}\mathcal{C}_{B}C_2\right]\varpi+\tilde{C}_F \le\frac{\varpi}{2}+\frac{\varpi}{2}=\varpi.\label{e648}
\end{multline}
From (\ref{e647}) and (\ref{e648}) we conclude that
$$\left\|\mathcal{H}_\epsilon\begin{pmatrix} \omega_0(\tau,m)\\ \omega_1(\tau,m)\end{pmatrix}\right\|\le\max\left\{\varpi,\varpi\right\}=\varpi.$$ 
The first statement is proved.

For the second part of the proof, let $\omega_{jk}\in \overline{B}(0,\varpi)\subseteq\hbox{Exp}^{q}_{(k,\beta,\mu,\alpha,\rho)}$, for $(j,k)\in\{0,1\}\times\{1,2\}$. We have 
$$\mathcal{H}_{\epsilon}\begin{pmatrix} \omega_{01}(\tau,m)\\\omega_{11}(\tau,m)\end{pmatrix}-\mathcal{H}_{\epsilon}\begin{pmatrix} \omega_{02}(\tau,m)\\\omega_{12}(\tau,m)\end{pmatrix}$$
equals
$$
\begin{pmatrix}\displaystyle
\mathcal{H}_{P}(\omega_{11}-\omega_{12})+\sum_{\ell=1}^{D-1}\epsilon^{\Delta_\ell-d_{\ell}}(\mathcal{H}_{\ell}(\omega_{01}-\omega_{02})+\delta_{\ell}\mathcal{H}_{\ell}(\omega_{11}-\omega_{12})) +\frac{1}{P_m(\tau)}\sum_{j=0}^{1}\tilde{b}_{j0}\star(\omega_{j1}-\omega_{j2})\\
\displaystyle 
\sum_{\ell=1}^{D-1}\epsilon^{\Delta_\ell-d_{\ell}}\mathcal{H}_{\ell}(\omega_{11}-\omega_{12})+\frac{1}{P_m(\tau)}\sum_{j=0}^{1}\tilde{b}_{j1}\star(\omega_{j1}-\omega_{j2})
\end{pmatrix}.
$$

An analogous reasoning as above yields
\begin{multline}
\left\|\mathcal{H}_{P}(\omega_{11}-\omega_{12})+\sum_{\ell=1}^{D-1}\epsilon^{\Delta_\ell-d_{\ell}}(\mathcal{H}_{\ell}(\omega_{01}-\omega_{02})+\delta_{\ell}\mathcal{H}_{\ell}(\omega_{11}-\omega_{12})) +\frac{1}{P_m(\tau)}\sum_{j=0}^{1}\tilde{b}_{j0}\star(\omega_{j1}-\omega_{j2})\right\|_{(k,\beta,\mu,\alpha,\rho)}\\
\le \frac{d_D}{k}\frac{1}{\mathfrak{D}_1}\max\{\frac{1}{C_D},\frac{1}{\mathfrak{D}_3}\}\left\|\omega_{11}-\omega_{12}\right\|_{(k,\beta,\mu,\alpha,\rho)}\\
\hfill+\sum_{\ell=1}^{D-1}\epsilon_0^{\Delta_\ell-d_{\ell}}C_{3,\ell}\left(\left\|\omega_{01}-\omega_{02}\right\|_{(k,\beta,\mu,\alpha,\rho)}+\delta_{\ell}\left\|\omega_{11}-\omega_{12}\right\|_{(k,\beta,\mu,\alpha,\rho)}\right)\\
\hfill+\sup_{m\in\R}\left(\frac{1}{|Q(im)|}\right)\frac{1}{C_D}\mathcal{C}_{B}C_2\left(\left\|\omega_{01}-\omega_{02}\right\|_{(k,\beta,\mu,\alpha,\rho)}+\left\|\omega_{11}-\omega_{12}\right\|_{(k,\beta,\mu,\alpha,\rho)}\right)\\
\le \left[\frac{d_D}{k}\frac{1}{\mathfrak{D}_1}\max\{\frac{1}{C_D},\frac{1}{\mathfrak{D}_3}\}+\sum_{\ell=1}^{D-1}\epsilon_0^{\Delta_\ell-d_{\ell}}C_{3,\ell}(1+\delta_{\ell})+\sup_{m\in\R}\left(\frac{1}{|Q(im)|}\right)\frac{2}{C_D}\mathcal{C}_{B}C_2\right]\\
\hfill\times\max_{j=0,1}(\left\|\omega_{j1}-\omega_{j2}\right\|_{(k,\beta,\mu,\alpha,\rho)})\\
\le\frac{1}{2}\max_{j=0,1}(\left\|\omega_{j1}-\omega_{j2}\right\|_{(k,\beta,\mu,\alpha,\rho)}),\label{e686}
\end{multline}
together with
\begin{multline}
\left\|\sum_{\ell=1}^{D-1}\epsilon^{\Delta_\ell-d_{\ell}}\mathcal{H}_{\ell}(\omega_{11}-\omega_{12})+\frac{1}{P_m(\tau)}\sum_{j=0}^{1}\tilde{b}_{j1}\star(\omega_{j1}-\omega_{j2})\right\|_{(k,\beta,\mu,\alpha,\rho)}\\
\le \sum_{\ell=1}^{D-1}\epsilon_0^{\Delta_\ell-d_{\ell}}C_{3,\ell}\left\|\omega_{11}-\omega_{12}\right\|_{(k,\beta,\mu,\alpha,\rho)}\hfill\\
+\sup_{m\in\R}\left(\frac{1}{|Q(im)|}\right)\frac{1}{C_D}\mathcal{C}_{B}C_2\left(\left\|\omega_{01}-\omega_{02}\right\|_{(k,\beta,\mu,\alpha,\rho)}+\left\|\omega_{11}-\omega_{12}\right\|_{(k,\beta,\mu,\alpha,\rho)}\right)\\
\le \left[\sum_{\ell=1}^{D-1}\epsilon_0^{\Delta_\ell-d_{\ell}}C_{3,\ell}+\sup_{m\in\R}\left(\frac{1}{|Q(im)|}\right)\frac{2}{C_D}\mathcal{C}_{B}C_2\right]\max_{j=0,1}(\left\|\omega_{j1}-\omega_{j2}\right\|_{(k,\beta,\mu,\alpha,\rho)})\\
\le\frac{1}{2}\max_{j=0,1}(\left\|\omega_{j1}-\omega_{j2}\right\|_{(k,\beta,\mu,\alpha,\rho)}).\label{e687}
\end{multline}
Regarding (\ref{e686}) and (\ref{e687}), we conclude that
\begin{multline*}
\left\|\mathcal{H}_{\epsilon}\begin{pmatrix} \omega_{01}(\tau,m)\\\omega_{11}(\tau,m)\end{pmatrix}-\mathcal{H}_{\epsilon}\begin{pmatrix} \omega_{02}(\tau,m)\\\omega_{12}(\tau,m)\end{pmatrix}\right\|\\
\le\max\{\frac{1}{2}\max_{j=0,1}\left\|\omega_{j1}-\omega_{j2}\right\|_{(k,\beta,\mu,\alpha,\rho)},\frac{1}{2}\max_{j=0,1}\left\|\omega_{j1}-\omega_{j2}\right\|_{(k,\beta,\mu,\alpha,\rho)}\}\\
=\frac{1}{2}\left\|\begin{pmatrix} \omega_{01}(\tau,m)\\\omega_{11}(\tau,m)\end{pmatrix}-\begin{pmatrix} \omega_{02}(\tau,m)\\\omega_{12}(\tau,m)\end{pmatrix}\right\|.
\end{multline*}
\end{proof}

\textbf{Remark:} We stress the availability of choice of the elements satisfying (\ref{e636}):
\begin{itemize}
\item Fix $Q, R_D$ to obtain $\mathfrak{D}_1$ and $\mathfrak{D}_2$ in Assumption (C). Fix all the parameters in the problem except $k$ and $d_D$.
\item Take $d_D$ to arrive at the bounds $\mathfrak{D}_3$ and $C_D$ provided in Corollary~\ref{coro495}. Then, fix $k$ and assume that the ratio $d_D/k$ is small enough so Assumption (D) holds (the first term in the sum of (\ref{e636}) can be as close to zero as needed). 
\item Choose small enough $\epsilon_0$. Assumption (A) guarantees that the terms of the sum in $\ell$ of (\ref{e636}) can be chosen as close to zero as needed.
\item Once the previous values are fixed, choose the coefficients $b_{jk}$ in such a way that $\mathcal{C}_B$ is small enough to have (\ref{e636}). 
\end{itemize}

\begin{example}
In order to illustrate the previous remark, we depart from a problem (\ref{e376}),(\ref{e377}) such that $\epsilon_0$ together with the constant $\mathcal{C}_{B}$ in (\ref{e321}), involved in the construction of the coefficients $b_{jk}$ for $j,k\in\{0,1\}$, are close to zero. Assume $d_D=1$, and take $Q(z)=1-z^2$, $R_{D}(z)=-2+z^2$. We observe that 
$$\mathfrak{D}_1=\frac{1}{2}\le \frac{Q(im)}{R_{D}(im)}=\frac{1+m^2}{2+m^2}\le 1=\mathfrak{D}_2.$$
Assume $d=0$.
The polynomial 
$$P_m(\tau)=m^2+1+(m^2+2)\tau$$
admits $q_0(m)=\frac{-m^2-1}{m^2+2}$ as its only root. We take $d=0$ and $S_d$ of small enough opening. $0<\rho<1/4$ and $\delta>0$ is large enough to fulfill the conditions at the beginning of Section~\ref{sec6}. Observe that $\hbox{arg}(q_0)=\pi$ for all $m\in\R$. Therefore, from the writing $\tau=ue^{i\theta}q_0(m)$ for all $\tau\in S_d$ one has that $\theta$ is close to $\pi$. The geometric conditions of Corollary~\ref{coro495} provide the following constants for $R=2$: 
$$\min_{u\in[0,2]}\frac{|1-ue^{i\pi}|}{1+u}\ge 1=:\mathfrak{D}_{31},$$
and
$$\min_{u\ge 2}\frac{u}{1+u}\frac{1}{2}=\frac{2}{3}\frac{1}{2}=:\mathfrak{D}_{32}.$$
We take $\mathfrak{D}_{3}=\min\{\mathfrak{D}_{31},\mathfrak{D}_{32}\}=\frac{1}{3}$. Observe that $C_D=1/2$ is a valid choice due to
$$|R_m(\tau)|=|m^2+1+(m^2+2)\tau|\ge m^2+1-\frac{1}{4}(m^2+2)=|Q(im)|\frac{3m^2+2}{4m^2+4}\ge |Q(im)|\frac{1}{2}.$$
In order that (\ref{e636}) holds, the valid values of $k$ are determined by
$$\frac{d_D}{k}\frac{1}{\mathfrak{D}_{1}}\max\{\frac{1}{C_D},\frac{1}{\mathfrak{D}_3}\}< \frac{1}{2},$$
i.e. $k>12$. A more accurate limit value for the valid values of $k$ can be given when providing more information on the geometric elements involved in the problem. From Assumption (A), the values of the other elements involved in the problem should satisfy in this situation that $\delta_{\ell}<1/12$ and $\Delta_{\ell}>d_{\ell}>1$ with $\hbox{deg}(R_{\ell})\le 2$ for all $1\le \ell\le D-1$. 
\end{example}

\section{Analytic solutions of the main problem}\label{sec7}

In this section we construct analytic solutions of the main problem under study. 

Let us consider the main equation (\ref{epral}), and assume the elements involved in the equation satisfy Assumptions (A), (B1), (B2) and (C). 

Let $\mathcal{E}\subseteq D(0,\epsilon_0)$ be a bounded sector with small opening. We also fix $\mathcal{T}\subseteq D(0,r_\mathcal{T})$ another bounded sector of small opening, for some $0<r_{\mathcal{T}}<1$. Both sectors have their vertex at the origin. The choice of $d\in\R$ is made according to the prescriptions in Section~\ref{sec6} in such a way that Assumption (D) is satisfied. Moreover, the freedom of choice for $d$, and consequently of $S_d$, allows us to assume that $d\neq\pi$ can also be chosen in such a way that for every $\epsilon\in\mathcal{E}$ and all $t\in\mathcal{T}$, then $\epsilon t\in\mathcal{R}_{d,\Delta}\cap D(0,r_1)$, for some small enough $\Delta,r_1>0$ (see Lemma~\ref{lema222} for the precise definition of this set). One may consider a subset of $\mathcal{R}_{d,\Delta}$ in such a way that it contains a sector not crossing the set of negative real numbers, provided that the opening of $\mathcal{E}$ and $\mathcal{T}$ are small enough.

The numbers $\epsilon_0,\varpi,\varsigma_b>0$ are chosen in accordance with Proposition~\ref{prop641}. We additionally choose $\alpha>0$ such that
\begin{equation}\label{e753}
\alpha<\frac{1}{2}-\frac{k}{\log(q)}\left(\log(1+2\delta/\rho)+\log(r_{\mathcal{T}})\right).
\end{equation}

The previous technical condition will be used in the proof of Theorem~\ref{teodif} at the time of giving upper bounds for the difference of two solutions in a common domain. Observe it is possible to choose such $\alpha$ if the radius $r_{\mathcal{T}}$ is small enough.

The first main result of the present work states the existence of an analytic solution of the main equation (\ref{epral}).

\begin{theo}\label{teo1}
The equation (\ref{epral}) admits a solution $u(t,z,\epsilon)$ of the form
\begin{equation}\label{e755}
u(t,z,\epsilon)=u_0(t,z,\epsilon)+u_1(t,z,\epsilon)\frac{\log(\epsilon t)}{\log(q)},
\end{equation}
holomorphic on $\mathcal{T}\times H_{\beta'}\times \mathcal{E}$, for every $0<\beta'<\beta$.
\end{theo}
\begin{proof}
For every $\epsilon\in\mathcal{E}$, Proposition~\ref{prop641} guarantees that $\mathcal{H}_\epsilon$ is a contractive map from a closed disc in a Banach space into itself. This entails the existence of a unique fixed point for $\mathcal{H}_{\epsilon}$. We observe this maps depends holomorphically on $\epsilon\in\mathcal{E}$. As a consequence, we arrive at the existence of $(\omega_0,\omega_1)\in(\overline{B}(0,\varpi))^2\subseteq(\hbox{Exp}_{(k,\beta,\mu,\alpha,\rho)}^{q})^2$ such that  $\mathcal{H}_\epsilon(\omega_0,\omega_1)=(\omega_0,\omega_1)$. Let us consider the functions 
$$(\tau,m,\epsilon)\mapsto \omega_{j}(\tau,m,\epsilon),$$
for $j\in\{0,1\}$, which are continuous on $(S_d\cup \overline{D}(0,\rho))\times\R\times\mathcal{E}$, holomorphic with respect to their first and third variable on $S_d\cup D(0,\rho)$ and $\mathcal{E}$, respectively. We observe that for $j\in\{0,1\}$
\begin{equation}\label{e767}
\sup_{\stackrel{\tau\in S_{d}\cup\overline{D}(0,\rho)}{m\in\R}}|\omega_j(\tau,m,\epsilon)|\le \varpi\frac{1}{(1+|m|)^{\mu}}e^{-\beta |m|}\exp\left(\frac{k}{2}\frac{\log^2|\tau+\delta|}{\log(q)}+\alpha\log|\tau+\delta|\right),
\end{equation}
for every $\tau\in S_d\cup\overline{D}(0,\rho)$, $m\in\R$ and $\epsilon\in\mathcal{E}$. Regarding the results on integral transforms in Section~\ref{secreview} one can define
$$u_j(t,z,\epsilon)=\frac{1}{(2\pi)^{1/2}}\frac{k}{\log(q)}\int_{-\infty}^{\infty}\int_{L_d}\frac{\omega_{j}(u,m,\epsilon)}{\Theta_{q^{1/k}}\left(\frac{u}{\epsilon t}\right)}\frac{du}{u}\exp(izm)dm,$$
for $j\in\{0,1\}$, which turns out to be a holomorphic function defined on $\mathcal{T}\times H_{\beta'}\times \mathcal{E}$, for any fixed $0<\beta'<\beta$. We observe from the strategy followed in the transformation of the main equation in Section~\ref{secproblemstrategy}, together with the properties of inverse Fourier and $q-$Laplace transform of order $k$ stated in Section~\ref{secreview}, that (\ref{e755}) is an analytic solution of (\ref{epral}).
\end{proof}

\section{Asymptotic study of the solutions}\label{sec9}

In this section we prove the existence of a formal solution to the main problem, and we state the asymptotic relation joining the analytic and the formal. This is possible by means of the application of a $q-$analog of the cohomological criteria known as Ramis-Sibuya Theorem. The classical result can be found in~\cite{hssi}, Lemma XI-2-6, whereas a $q-$analog of this result was recalled in Section~\ref{sec23}.

We depart from the main equation (\ref{epral}), with its elements satisfying Assumptions (A), (B2) and (C). The forcing term, $f$ is constructed in the same fashion as in Assumption (B1), with direction $d$ be chosen among the elements in the finite set $(d_p)_{0\le p\le\zeta-1}$ to be determined.

We first consider an appropriate geometric framework for the domain of definition of the perturbation parameter $\epsilon$.

\begin{defin}
Let $\zeta\ge2$ be an integer and $(\mathcal{E}_p)_{0\le p\le \zeta-1}$ be a good covering in $\C^{\star}$ (see Definition~\ref{defigc}) and an open bounded sector $\mathcal{T}$ with vertex at the origin and radius $r_{\mathcal{T}}>0$. We consider a family of unbounded sectors $(S_{d_{p}})_{0\le p\le\zeta-1}$, with vertex at the origin and bisecting direction $d_{p}\in\R$ for every $0\le p\le\zeta-1$. The direction $d_p$ is chosen so that the geometric conditions for $d$ in Section~\ref{sec6} are satisfied. More precisely, for all $0\le p\le \zeta-1$, one assumes the following statements regarding $d_p$:
\begin{itemize}
\item[(iii)] $S_{d_{p}}$ avoids the roots of $P_m$ defined in (\ref{e434}) and $\hbox{dist}(S_{d_p}\cup D(0,\rho),-\delta)\ge1$.
\item[(iv)] Assumption (D) holds.
\item[(v)] For every $\epsilon\in\mathcal{E}_p$ and all $t\in\mathcal{T}$, $\epsilon t\in\mathcal{R}_{d_p,\Delta}\cap D(0,r_1)$ for some small enough $\Delta,r_1>0$.
\end{itemize}
The family $\{(\mathcal{R}_{d_p,\Delta})_{0\le p\le\zeta-1},D(0,\rho),\mathcal{T}\}$ is said to be associated to the good covering $(\mathcal{E}_p)_{0\le p\le \zeta-1}$. 
\end{defin}

Let $(\mathcal{E}_{p})_{0\le p\le \zeta-1}$ be a good covering in $\C^{\star}$ and fix a set $\{(\mathcal{R}_{d_p,\Delta})_{0\le p\le\zeta-1},D(0,\rho),\mathcal{T}\}$, associated to the previous good covering. 

We choose the numbers $\epsilon_0,\varpi,\varsigma_b>0$ in order that Proposition~\ref{prop641} holds for every choice of $\mathcal{E}_p$ among the elements in the good covering. Due to the nature of the forcing term $f$, for every $0\le p\le \zeta-1$, we may vary the value of $d$ in Assumption (B1) among the elements in $(d_p)_{0\le p\le \zeta-1}$.  We recall that any choice of $d_p\in\R$ can be made which does not vary the definition of $f$ due to its polynomial nature.

Fix $0<\beta'<\beta$ and for every $0\le p\le \zeta-1$ we build the solution of (\ref{epral}) in the form 
\begin{equation}\label{e859}
u_p(t,z,\epsilon)=u_{0,p}(t,z,\epsilon)+u_{1,p}(t,z,\epsilon)\frac{\log(\epsilon t)}{\log(q)},
\end{equation}
where $u_{0,p}$ and $u_{1,p}$ are constructed following Theorem~\ref{teo1}. More precisely one has that
$$u_{j,p}(t,z,\epsilon)=\frac{1}{(2\pi)^{1/2}}\frac{k}{\log(q)}\int_{-\infty}^{\infty}\int_{L_{d_p}}\frac{\omega_{j,p}(u,m,\epsilon)}{\Theta_{q^{1/k}}\left(\frac{u}{\epsilon t}\right)}\frac{du}{u}\exp(izm)dm,$$
for $j\in\{0,1\}$, and where $\omega_{j,p}(u,m,\epsilon)$ is constructed from the fixed point argument, mimicking Section~\ref{sec6}. Therefore, it satisfies that 
$$\sup_{\stackrel{\tau\in S_{d_p}\cup\overline{D}(0,\rho)}{m\in\R}}|\omega_{j,p}(\tau,m,\epsilon)|\le \varpi\frac{1}{(1+|m|)^{\mu}}e^{-\beta |m|}\exp\left(\frac{k}{2}\frac{\log^2|\tau+\delta|}{\log(q)}+\alpha\log|\tau+\delta|\right),$$
for every $\epsilon\in\mathcal{E}_p\cap\mathcal{E}_{p+1}$. We recall that $u_p\in\mathcal{O}(\mathcal{T}\times H_{\beta'}\times \mathcal{E}_p)$ for every $0\le p\le \zeta-1$.

\begin{theo}\label{teodif}
In the previous situation, there exist constants $\tilde{K}>0$ and $\tilde{\alpha}\in\R$ such that 
$$\sup_{t\in\mathcal{T},z\in H_{\beta'}}|u_{j,p+1}(t,z,\epsilon)-u_{j,p}(t,z,\epsilon)|\le \tilde{K}\exp\left(-\frac{k}{2\log(q)}\log^2|\epsilon|\right)|\epsilon|^{\tilde{\alpha}} $$
for $\epsilon\in \mathcal{E}_{p}\cap\mathcal{E}_{p+1}$, and all $0\le p\le \zeta-1$ (by identifying $u_{j,\zeta}$ with $u_{j,0}$), for $j\in\{0,1\}$. 
\end{theo}
\begin{proof}
Fix $0\le p\le \zeta-1$. We observe that for $j\in\{0,1\}$, every $t\in\mathcal{T}$, $z\in H_{\beta'}$ and $\epsilon\in \mathcal{E}_{p}\cap\mathcal{E}_{p+1}$ the difference $u_{j,p+1}(t,z,\epsilon)- u_{j,p}(t,z,\epsilon)$ is given by
\begin{equation}\label{e818}
\frac{1}{(2\pi)^{1/2}}\frac{k}{\log(q)}\int_{-\infty}^{\infty}\left[\int_{L_{d_{p+1}}}\frac{\omega_{j,p+1}(u,m,\epsilon)}{\Theta_{q^{1/k}}\left(\frac{u}{\epsilon t}\right)}\frac{du}{u}-\int_{L_{d_p}}\frac{\omega_{j,p}(u,m,\epsilon)}{\Theta_{q^{1/k}}\left(\frac{u}{\epsilon t}\right)}\frac{du}{u}\right]\exp(izm)dm.
\end{equation}

An analogous reasoning as in the construction of the analytic solutions to the main problem, substituting the Banach space $\hbox{Exp}_{(k,\beta,\mu,\alpha,\rho)}^{q}$ by the same Banach space substituting the bounds on the set $\overline{D}(0,\rho)$ with respect to $\tau$, instead of $\overline{D}(0,\rho)\cup S_d$ guarantees the existence of a unique solution of the main problem in such Banach space. On the other hand, the restriction of $\omega_{j,p}$ and $\omega_{j,p+1}$ to $\overline{D}(0,\rho)$ are both solutions to the same problem. Therefore, both functions coincide in the set $\overline{D}(0,\rho)\times\R\times (\mathcal{E}_{p}\cap\mathcal{E}_{p+1})$. Let us denote $\omega_{j,p,p+1}$ the common function defined in $\overline{D}(0,\rho)\times\R\times (\mathcal{E}_{p}\cap\mathcal{E}_{p+1})$. 

A path deformation of the Laplace integrals defining the difference in (\ref{e818}) together with Cauchy theorem allows to write this difference in the form
\begin{multline}
\frac{1}{(2\pi)^{1/2}}\frac{k}{\log(q)}\int_{-\infty}^{\infty}\left[\int_{L_{d_{p+1},\rho/2}}\frac{\omega_{j,p+1}(u,m,\epsilon)}{\Theta_{q^{1/k}}\left(\frac{u}{\epsilon t}\right)}\frac{du}{u}+\int_{C_{j,p,p+1}(\rho/2)}\frac{\omega_{j,p,p+1}(u,m,\epsilon)}{\Theta_{q^{1/k}}\left(\frac{u}{\epsilon t}\right)}\frac{du}{u}\right.\\
\left.-\int_{L_{d_p},\rho/2}\frac{\omega_{j,p}(u,m,\epsilon)}{\Theta_{q^{1/k}}\left(\frac{u}{\epsilon t}\right)}\frac{du}{u}\right]\exp(izm)dm=:I_1+I_2-I_3,\label{e836}
\end{multline}
where the paths are given by $L_{d_{h},\rho/2}=\{re^{id_h}: r\in[\rho/2,\infty)\}$ for $h\in\{p,p+1\}$, and the arc $C_{j,p,p+1}(\rho/2)=\{\rho/2e^{i\theta}:\theta\in[d_p,d_{p+1}]\}$.
We now provide upper estimates for $|I_\ell|$ for $\ell=1,2,3$.

First, observe from (\ref{e767}) and Lemma~\ref{lema194} that 
\begin{multline*}
|I_1|=\left|\frac{1}{(2\pi)^{1/2}}\frac{k}{\log(q)}\int_{-\infty}^{\infty}\int_{L_{d_{p+1},\rho/2}}\frac{\omega_{j,p+1}(u,m,\epsilon)}{\Theta_{q^{1/k}}\left(\frac{u}{\epsilon t}\right)}\frac{du}{u}e^{izm}dm\right|\\
\le\frac{1}{(2\pi)^{1/2}}\frac{k}{\log(q)}\int_{-\infty}^{\infty}\int_{\rho/2}^{\infty}\frac{|\omega_{j,p+1}(re^{id_{p+1}},m,\epsilon)|}{\left|\Theta_{q^{1/k}}\left(\frac{re^{id_{p+1}}}{\epsilon t}\right)\right|}\frac{dr}{r}e^{-m\hbox{Im}(z)}dm\\
\le \frac{\varpi}{(2\pi)^{1/2}}\frac{k}{\log(q)}\frac{1}{C_{q,k}\Delta}\left[\int_{-\infty}^{\infty}\frac{1}{(1+|m|)^{\mu}}e^{-|m|(\beta- \beta')}dm\right]\\
\times|\epsilon t|^{1/2}\int_{\rho/2}^{\infty} \frac{\exp\left(\frac{k}{2}\frac{\log^2|re^{id_{p+1}}+\delta|}{\log(q)}+\alpha\log|re^{id_{p+1}}+\delta|\right)}{\exp\left(\frac{k}{2}\frac{\log^2(r/|\epsilon t|)}{\log(q)}\right)}\frac{dr}{r^{3/2}}.
\end{multline*}
The previous to the last line is a constant. We proceed to upper estimate the last line of the previous expression by observing that  
\begin{multline}
|\epsilon t|^{1/2}\int_{\rho/2}^{\infty} \frac{\exp\left(\frac{k}{2}\frac{\log^2|re^{id_{p+1}}+\delta|}{\log(q)}+\alpha\log|re^{id_{p+1}}+\delta|\right)}{\exp\left(\frac{k}{2}\frac{\log^2(r/|\epsilon t|)}{\log(q)}\right)}\frac{dr}{r^{3/2}}\\
\le |\epsilon t|^{1/2}\int_{\rho/2}^{\infty} \exp\left(\frac{k}{2\log(q)}\left(\log^2(r+\delta)-\log^2\left(\frac{r}{|\epsilon t|}\right)\right)\right)\frac{(r+\delta)^{\alpha}}{r^{3/2}}dr\\
=\left(|t|^{1/2}e^{-\frac{k}{2\log(q)}\log^2|t|}\right)|\epsilon|^{1/2}e^{\frac{k\log^2(1+2\delta/\rho)}{2\log(q)}}\\
\times\int_{\rho/2}^{\infty} \exp\left(\frac{k}{2\log(q)}\left(-2\log|\epsilon|\log|t|+2\log(r)\log|\epsilon|+2\log(r)\log|t|\right)\right)\frac{(r+\delta)^{\alpha+\frac{k\log(1+2\delta/\rho)}{\log(q)}}}{r^{3/2}}dr\\
\times\exp\left(-\frac{k}{2\log(q)}\log^2|\epsilon|\right).\label{e847}
\end{multline}

The choice $0<\epsilon_0<1$ and $0<r_{\mathcal{T}}<1$ yields 
$$\exp\left(-\frac{k}{\log(q)}\log|\epsilon|\log|t|\right)\le|\epsilon|^{-\frac{k}{\log(q)}\log(r_{\mathcal{T}})},$$ 
$$\exp\left(\frac{k}{\log(q)}\log(r)\log|\epsilon|\right)\le|\epsilon|^{\frac{k}{\log(q)}\log(\rho/2)},$$ 
for every $t\in\mathcal{T}$, $\epsilon\in\mathcal{E}_p\cap\mathcal{E}_{p+1}$ and $r\ge \rho/2$, together with
$$\exp\left(\frac{k}{\log(q)}\log(r)\log|t|\right)\le|t|^{\frac{k}{\log(q)}\log(\rho/2)},\quad \rho/2\le r\le 1,$$
and 
$$\exp\left(\frac{k}{\log(q)}\log(r)\log|t|\right)\le r^{\frac{k}{\log(q)}\log(r_{\mathcal{T}})},\quad r\ge1,$$
for all $t\in\mathcal{T}$.
This entails that (\ref{e847}) is upper bounded by
\begin{multline*}\left(|t|^{1/2}e^{-\frac{k}{2\log(q)}\log^2|t|}\right)|\epsilon|^{1/2+\frac{k}{\log(q)}(-\log(r_{\mathcal{T}})+\log(\rho/2))}e^{\frac{k\log^2(1+2\delta/\rho)}{2\log(q)}}\\
\times\left[\int_{\rho/2}^{1}|t|^{\frac{k}{\log(q)}\log(\rho/2)} \frac{(r+\delta)^{\alpha+\frac{k\log(1+2\delta/\rho)}{\log(q)}}}{r^{3/2}}dr+\int_{1}^{\infty}r^{\frac{k}{\log(q)}\log(r_\mathcal{T})} \frac{(r+\delta)^{\alpha+\frac{k\log(1+2\delta/\rho)}{\log(q)}}}{r^{3/2}}dr\right]\\
\times\exp\left(-\frac{k}{2\log(q)}\log^2|\epsilon|\right).
\end{multline*}

It is straight to check the existence of a positive constant $\tilde{K}_1$ such that
\begin{equation}\label{e872}
\sup_{|t|>0}\max\{|t|^{1/2+\frac{k}{\log(q)}\log(\rho/2)},|t|^{1/2}\}e^{-\frac{k}{2\log(q)}\log^2|t|}\le \tilde{K}_1.
\end{equation}
In addition to this, there exists a constant $\tilde{K}_2>0$ such that
$$\int_{\rho/2}^{1} \frac{(r+\delta)^{\alpha+\frac{k\log(1+2\delta/\rho)}{\log(q)}}}{r^{3/2}}dr\le \tilde{K}_2.$$
and from the choice made on the value of $\alpha$ in (\ref{e753}), there exists another constant $\tilde{K}_3>0$ such that
$$\int_{1}^{\infty}r^{\frac{k}{\log(q)}\log(r_\mathcal{T})} \frac{(r+\delta)^{\alpha+\frac{k\log(1+2\delta/\rho)}{\log(q)}}}{r^{3/2}}dr\le \tilde{K}_3.$$
We take $\tilde{K}_{I_1}:=\tilde{K}_1e^{\frac{k\log^2(1+2\delta/\rho)}{2\log(q)}}(\tilde{K}_2+\tilde{K}_3)$ and $\tilde{\alpha}_{I_1}:=\frac{1}{2}+\frac{k}{\log(q)}(\log(\rho/2)-\log(r_{\mathcal{T}}))$ to conclude that
\begin{equation}\label{e001}
|I_1|\le \tilde{K}_{I_1} \exp\left(-\frac{k}{2\log(q)}\log^2|\epsilon|\right) |\epsilon|^{\tilde{\alpha}_{I_1}}.
\end{equation}
An analogous reasoning guarantees that
\begin{equation}\label{e003}
|I_3|\le \tilde{K}_{I_3} \exp\left(-\frac{k}{2\log(q)}\log^2|\epsilon|\right) |\epsilon|^{\tilde{\alpha}_{I_3}},
\end{equation}
for $\tilde{K}_{I_3}:=\tilde{K}_{I_1}$ and $\tilde{\alpha}_{I_3}:=\tilde{\alpha}_{I_1}$.

We conclude the proof by estimating $|I_2|$. Analogous computations as above yield
\begin{multline*}
|I_2|=\left|\frac{1}{(2\pi)^{1/2}}\frac{k}{\log(q)}\int_{-\infty}^{\infty}\int_{C_{j,p,p+1}(\rho/2)}\frac{\omega_{j,p,p+1}(u,m,\epsilon)}{\Theta_{q^{1/k}}\left(\frac{u}{\epsilon t}\right)}\frac{du}{u}\exp(izm)dm\right|\\
\le\frac{1}{(2\pi)^{1/2}}\frac{k}{\log(q)}\int_{-\infty}^{\infty}\int_{d_p}^{d_{p+1}}\frac{|\omega_{j,p,p+1}(\rho/2e^{i\theta},m,\epsilon)|}{\left|\Theta_{q^{1/k}}\left(\frac{\rho/2e^{i\theta}}{\epsilon t}\right)\right|}d\theta e^{-m\hbox{Im}(z)}dm\\
\le\frac{\varpi}{(2\pi)^{1/2}}\frac{k}{\log(q)}\frac{1}{C_{q,k}\Delta}\left[\int_{-\infty}^{\infty}\frac{1}{(1+|m|)^{\mu}}e^{-|m|(\beta- \beta')}dm\right]\\
\times|\epsilon t|^{1/2}(2/\rho)^{1/2}\int_{d_p}^{d_{p+1}} \frac{\exp\left(\frac{k}{2}\frac{\log^2|\rho/2e^{i\theta}+\delta|}{\log(q)}+\alpha\log|\rho/2e^{i\theta}+\delta|\right)}{\exp\left(\frac{k}{2}\frac{\log^2(\rho/(2|\epsilon t|))}{\log(q)}\right)}d\theta\\
\le\frac{\varpi}{(2\pi)^{1/2}}\frac{k}{\log(q)}\frac{1}{C_{q,k}\Delta}\left[\int_{-\infty}^{\infty}\frac{1}{(1+|m|)^{\mu}}e^{-|m|(\beta- \beta')}dm\right]|d_{p+1}-d_p|\\
\exp\left(\frac{k}{2}\frac{\log^2(\rho/2+\delta)}{\log(q)}+\alpha\log(\rho/2+\delta)\right)\times|\epsilon|^{1/2}r_{\mathcal{T}}^{1/2}(2/\rho)^{1/2} \exp\left(-\frac{k}{2}\frac{\log^2(\rho/(2|\epsilon t|))}{\log(q)}\right).
\end{multline*}
At this point, we apply (\ref{e872}) and the fact that $0<\epsilon_0,r_\mathcal{T}<1$ to arrive at
\begin{multline*}
|\epsilon|^{1/2}\exp\left(-\frac{k}{2}\frac{\log^2(\rho/(2|\epsilon t|))}{\log(q)}\right)\\
=\exp\left(-\frac{k}{2\log(q)}\log^2(\rho/2)\right)|\epsilon|^{\frac{1}{2}+\frac{k}{\log(q)}\log(\rho/2)}\left[|t|^{\frac{k}{\log(q)}\log(\rho/2)}\exp\left(-\frac{k}{2\log(q)}\log^2|t|\right)\right]\\
\times\exp\left(\frac{k}{2\log(q)}(-\log^2|\epsilon|-2\log|\epsilon|\log|t|)\right)\\
\le \exp\left(-\frac{k}{2\log(q)}\log^2(\rho/2)\right)|\epsilon|^{\frac{1}{2}+\frac{k}{\log(q)}\log(\rho/2)}\tilde{K}_1
\exp\left(\frac{k}{2\log(q)}(-\log^2|\epsilon|)\right),
\end{multline*}
to obtain the existence of $\tilde{K}_{I_2}>0$ and  $\tilde{\alpha}_{I_2}\in\R$ with

\begin{equation}\label{e002}
|I_2|\le \tilde{K}_{I_2} \exp\left(-\frac{k}{2\log(q)}\log^2|\epsilon|\right) |\epsilon|^{\tilde{\alpha}_{I_2}}.
\end{equation}

In accordance to (\ref{e836}) and the estimates (\ref{e003}), (\ref{e001}) and (\ref{e002}), the result follows.

\end{proof}

We are in conditions to state the main result of the present work.

Let $0<\beta'<\beta$ and consider be the Banach space $\mathbb{E}$ of bounded holomorphic functions on $\mathcal{T}\times H_{\beta'}$, with the norm of the supremum.

\begin{theo}\label{teopral}
There exist $\hat{u}_0,\hat{u}_1\in\mathbb{E}[[\epsilon]]$ such that the function $\epsilon\mapsto u_{j,p}$, defined in (\ref{e859}), admits $\hat{u}_j$ as its $q-$Gevrey asymptotic expansion of order $1/k$ on $\mathcal{E}_p$ for every $0\le p\le \zeta-1$ and $j\in\{0,1\}$.

In addition to this, the formal expression
\begin{equation}\label{e974}
\hat{u}(t,z,\epsilon)=\hat{u}_0(t,z,\epsilon)+\hat{u}_1(t,z,\epsilon)\frac{\log(\epsilon t)}{\log(q)}
\end{equation}
is a formal solution of
\begin{multline}
Q(\partial_z)\hat{u}(t,z,\epsilon)=(\epsilon t)^{d_D}\sigma_{q;t}^{\frac{d_D}{k}}R_D(\partial_{z})\hat{u}(t,z,\epsilon)+\sum_{\ell=1}^{D-1}\epsilon^{\Delta_\ell}t^{d_{\ell}}c_{\ell}(z,\epsilon)R_{\ell}(\partial_z)\sigma_{q;t}^{\delta_{\ell}}\hat{u}(t,z,\epsilon)\\
+f(t,z,\epsilon)+\left(b_{00}(z,\epsilon) +b_{01}(z,\epsilon)\frac{\log(\epsilon t)}{\log(q)}\right)\left[\hat{u}(t,z,\epsilon)-\left(\frac{1}{2\pi i}(\gamma_{\epsilon}^{*}-\hbox{id})\hat{u}(t,z,\epsilon)\right)\log(\epsilon t)\right]\\
+\left(b_{10}(z,\epsilon)+b_{11}(z,\epsilon)\frac{\log(\epsilon t)}{\log(q)}\right)\left[\frac{\log(q)}{2\pi i}(\gamma_{\epsilon}^{*}-\hbox{id})\hat{u}(t,z,\epsilon)\right]\label{epralformal}.
\end{multline}
with $\gamma_{\epsilon}^{*}$ being the formal monodromy operator defined in (\ref{e125}).
\end{theo}

\begin{proof}
Regarding Theorem~\ref{teodif} and Theorem ($q$-RS) in Section~\ref{sec23} there exist $\hat{u}_{0},\hat{u}_1\in\mathbb{E}[[\epsilon]]$ such that $u_{j,p}$ admits $\hat{u}_j$ as its $q-$Gevrey asymptotic expansion of order $1/k$ in $\mathcal{E}_p$, for every $0\le p\le \zeta-1$ and $j\in\{0,1\}$. 

It only rests to prove that (\ref{e974}) is a formal solution to (\ref{epralformal}). Let us write
\begin{equation}\label{e1010}
\hat{u}_j(t,z,\epsilon):=\sum_{n\ge0}\tilde{u}_{j,n}(t,z)\frac{\epsilon^n}{n!},\quad j\in\{0,1\},
\end{equation}
with $\tilde{u}_{j,n}\in\mathbb{E}$.  In view of Proposition~\ref{prop288} one has that for every $0\le p\le \zeta-1$ and every $\tilde{\mathcal{E}}_{p}\prec\mathcal{E}_p$
\begin{equation}\label{e992}
\lim_{\epsilon\to0,\epsilon\in\tilde{\mathcal{E}}_p}\left\|\partial_\epsilon^mu_{j,p}(t,z,\epsilon)-\tilde{u}_{j,m}(t,z)\right\|_{\mathbb{E}}=0,
\end{equation} 
for $j\in\{0,1\}$ and every integer $m\ge0$. The fact that $u_p$ is a solution of (\ref{epral}) guarantees that $u_{0,p}$ is a solution of (\ref{e235}) defined in $\mathcal{T}\times H_{\beta'}\times\mathcal{E}_p$ and  $u_{1,p}$ is a solution of (\ref{e236}) in $\mathcal{T}\times H_{\beta'}\times\mathcal{E}_p$. 

Let $m\ge\max\{d_{D},\max_{\ell=1,\ldots,D-1}\{\Delta_{\ell}\}\}$. Let $\tilde{\mathcal{E}}_p\prec\mathcal{E}_p$. We plug $u_0$ into (\ref{e235}) and apply the operator $\partial_\epsilon^{m}$ at both sides of the equation to arrive at
\begin{multline}
Q(\partial_z)\partial_{\epsilon}^{m}u_{0,p}(t,z,\epsilon)=\sum_{m_1+m_2=m}\frac{m!}{m_1!m_2!}(\partial_{\epsilon}^{m_1}\epsilon^{d_D}) t^{d_D}R_{D}(\partial_z)\\
\hfill\times\left[\partial_{\epsilon}^{m_2}u_{0,p}(q^{\frac{d_D}{k}}t,z,\epsilon)+\frac{d_D}{k}\partial_{\epsilon}^{m_2}u_{1,p}(q^{\frac{d_D}{k}}t,z,\epsilon)\right]\\
+\sum_{\ell=1}^{D-1}\sum_{m_1+m_2+m_3=m}\frac{m!}{m_1!m_2!m_3!}(\partial_{\epsilon}^{m_1}\epsilon^{\Delta_\ell})t^{d_{\ell}}(\partial_{\epsilon}^{m_2}c_{\ell}(z,\epsilon))R_{\ell}(\partial_z)\\
\hfill\times\left[\partial_{\epsilon}^{m_3}u_{0,p}(q^{\delta_{\ell}}t,z,\epsilon)+\delta_{\ell}\partial_{\epsilon}^{m_3}u_{1,p}(q^{\delta_{\ell}}t,z,\epsilon)\right]\\
+\partial_{\epsilon}^{m}f_0(t,z,\epsilon)+\sum_{j=0}^{1}\sum_{m_1+m_2=m}\frac{m!}{m_1!m_2!}(\partial_{\epsilon}^{m_1}b_{j0}(z,\epsilon))(\partial_{\epsilon}^{m_2}u_{j,p}(t,z,\epsilon)),\label{e1021}
\end{multline}
for all $(t,z,\epsilon)\in\mathcal{T}\times H_{\beta'}\times\mathcal{E}_p$. 
Let us write $c_{\ell}(z,\epsilon)=\sum_{n\ge0}\frac{c_{\ell,n}(z)}{n!}\epsilon^n$
for $1\le \ell\le D-1$ with $c_{\ell,n}\in E_{(\beta,\mu)}$. Also, we write $b_{jk}(z,\epsilon)=\sum_{n\ge0}\frac{b_{jk,n}(z)}{n!}\epsilon^n$ for $j,k\in\{0,1\}$, with $b_{jk,n}\in E_{(\beta,\mu)}$ for all $n\ge0$. Observe from the holomorphy of $c_{\ell}$ and $b_{jk}$ at the origin that 
$$\lim_{\epsilon\to0}\partial_\epsilon^{m_1}c_{\ell}(z,\epsilon)=c_{\ell,m_1}(z),$$
and
$$\lim_{\epsilon\to0}\partial_\epsilon^{m_1}b_{jk}(z,\epsilon)=b_{jk,m_1}(z)$$
for all nonnegative integer $m_1$.
We also write 
$$f_{h}(t,z,\epsilon)=\sum_{n\ge0}f_{h,n}(t,z)\frac{\epsilon^n}{n!},$$
for $h\in\{0,1\}$ and $(t,z)\in \C\times H_{\beta'}$, arriving at
$$\lim_{\epsilon\to0}\partial_\epsilon^{m_1}f_h(t,z,\epsilon)=f_{h,m_1}(t,z)$$
for all nonnegative integer $m_1$.

Taking the previous facts into account and (\ref{e992}), we get that for all $\epsilon\in\tilde{\mathcal{E}}_p$ one can tend $\epsilon\to0$ in (\ref{e1021}) to arrive at
\begin{multline}
Q(\partial_z)\frac{\tilde{u}_{0,m}(t,z)}{m!}= t^{d_D}R_{D}(\partial_z)\left[\frac{\tilde{u}_{0,m-d_D}(q^{\frac{d_D}{k}}t,z)}{(m-d_D)!}+\frac{d_D}{k}\frac{\tilde{u}_{1,m-d_D}(q^{\frac{d_D}{k}}t,z)}{(m-d_D)!}\right]\\
+\sum_{\ell=1}^{D-1}\sum_{m_2+m_3=m-\Delta_\ell}t^{d_{\ell}}\frac{c_{\ell,m_2}(z)}{m_2!}R_{\ell}(\partial_z)\left[\frac{\tilde{u}_{0,m_3}(q^{\delta_{\ell}}t,z)}{m_3!}+\delta_{\ell}\frac{\tilde{u}_{1,m_3}(q^{\delta_{\ell}}t,z)}{m_3!}\right]\\
+\frac{f_{0,m}(t,z)}{m!}+\sum_{j=0}^{1}\sum_{m_1+m_2=m}\frac{b_{j0,m_1}(z)}{m_1!}\frac{\tilde{u}_{j,m_2}(t,z)}{m_2!}.\label{e1021b}
\end{multline}

In an analogous manner we derive from (\ref{e236}) the equation

\begin{multline*}
Q(\partial_z)\partial_\epsilon^mu_{1,p}(t,z,\epsilon)=\sum_{m_1+m_2=m}\frac{m!}{m_1!m_2!}(\partial_\epsilon^{m_1}\epsilon^{d_D})t^{d_D}R_{D}(\partial_z)\partial_\epsilon^{m_2}u_{1,p}(q^{\frac{d_D}{k}}t,z,\epsilon)\\
+\sum_{\ell=1}^{D-1}\sum_{m_1+m_2+m_3=m}\frac{m!}{m_1!m_2!m_3!}(\partial_\epsilon^{m_1}\epsilon^{\Delta_\ell})t^{d_{\ell}}(\partial_\epsilon^{m_2}c_{\ell}(z,\epsilon))R_{\ell}(\partial_z)\partial_\epsilon^{m_3}u_1(q^{\delta_{\ell}}t,z,\epsilon)\\
+\partial_\epsilon^mf_1(t,z,\epsilon)+\sum_{j=0}^{1}\sum_{m_1+m_2=m}\frac{m!}{m_1!m_2!}(\partial_\epsilon^{m_1}b_{j1}(z,\epsilon))(\partial_\epsilon^{m_2}u_j(t,z,\epsilon)),
\end{multline*}

which leads to 

\begin{multline}
Q(\partial_z)\frac{\tilde{u}_{1,m}(t,z)}{m!}=t^{d_D}R_{D}(\partial_z)\frac{\tilde{u}_{1,m-d_D}(q^{\frac{d_D}{k}}t,z)}{(m-d_D)!}\\
+\sum_{\ell=1}^{D-1}\sum_{m_2+m_3=m-\Delta_\ell}t^{d_{\ell}}\frac{c_{\ell,m_2}(z)}{m_2!}R_{\ell}(\partial_z)\frac{\tilde{u}_{1,m_3}(q^{\delta_{\ell}}t,z)}{m_3!}\\
+\frac{f_{1,m}(t,z)}{m!}+\sum_{j=0}^{1}\sum_{m_1+m_2=m}\frac{b_{j1,m_1}(z)}{m_1!}\frac{\tilde{u}_{j,m_2}(t,z)}{m_2!}.\label{e1021c}
\end{multline}

On the other hand, it is straight to check that the formal power series defined by (\ref{e1010}) formally satisfy (\ref{e235}) and (\ref{e236}) by direct inspection. Indeed, one arrives at (\ref{e1021b}) and (\ref{e1021c}), respectively. This allows us to conclude that (\ref{e974}) is a formal solution of (\ref{epral}).
\end{proof}

\section{A problem associated to a triangular coupled system}\label{sec10}

In this section, we consider a simplified version of the main equation under study, for which the coupled system appearing in the construction of the analytic solution turns out to be triangularized, and the solution can be found in a more simpler way. Indeed, the assumptions to be made can be relaxed as we may observe.

The problem is presented, focusing on that parts in which there is a remarkable change in the procedure or in the results. Otherwise, the steps are skipped.

As in the general case, we fix an integer number $D\ge 2$ and $q\in\R$ with $q>1$. 

Let $d_{\ell},\Delta_{\ell}$  for $1\le \ell\le D-1$ and $d_D$ be non negative integers, and $\delta_{\ell}$ a non negative rational number for $1\le \ell\le D-1$. Let $\beta>0$ and fix $0<\beta'<\beta$. We also choose $0<\epsilon_0<1$, and $\mu>1$.  These parameters satisfy Assumption (A). We also take a small $\epsilon_0>0$. The polynomials $R_{\ell}$ for $1\le \ell\le D$ and $Q$ satisfy Assumptions (A) and (C). 

The forcing term $f(t,z,\epsilon)$ and the coefficients $c_{\ell}(z,\epsilon)$ for $1\le \ell\le D-1$ and $b_{jk}(z,\epsilon)$ for $j,k\in\{0,1\}$ with $(j,k)\neq(0,1)$ are determined by Assumptions (B1) and (B2). Observe that the only difference with respect to the general case is that we assume in this section that
$$b_{01}\equiv 0,$$
it is to say, there is an absence of logarithmic term at one of the coefficients of the operators involving the formal monodromy. The simplified version of the equation under study is 
\begin{multline}
Q(\partial_z)u(t,z,\epsilon)=(\epsilon t)^{d_D}\sigma_{q;t}^{\frac{d_D}{k}}R_D(\partial_{z})u(t,z,\epsilon)+\sum_{\ell=1}^{D-1}\epsilon^{\Delta_\ell}t^{d_{\ell}}c_{\ell}(z,\epsilon)R_{\ell}(\partial_z)\sigma_{q;t}^{\delta_{\ell}}u(t,z,\epsilon)\\
+f(t,z,\epsilon)+b_{00}(z,\epsilon)\left[u(t,z,\epsilon)-\left(\frac{1}{2\pi i}(\gamma_{\epsilon}^{*}-\hbox{id})u(t,z,\epsilon)\right)\log(\epsilon t)\right]\\
+\left(b_{10}(z,\epsilon)+b_{11}(z,\epsilon)\frac{\log(\epsilon t)}{\log(q)}\right)\left[\frac{\log(q)}{2\pi i}(\gamma_{\epsilon}^{*}-\hbox{id})u(t,z,\epsilon)\right]\label{eprals}.
\end{multline}

Assumption (B2) is substituted by the following assumption.

\vspace{0.3cm}

\textbf{Assumption (B2'):} The coefficients $c_\ell(z,\epsilon)$ for $1\le \ell\le D-1$, and $b_{jk}(z,\epsilon)$ for $j,k\in\{0,1\}$ with $(j,k)\neq (0,1)$ are holomorphic functions defined in $H_{\beta'}\times D(0,\epsilon_0)$ which belong to $E_{(\beta,\mu)}$, and are constructed as the inverse Fourier transform of functions $m\mapsto C_{\ell}(m,\epsilon)$ and $m\mapsto \tilde{b}_{jk}(m,\epsilon)$ respectively, holomorphic with respect to $\epsilon\in D(0,\epsilon_0)$ such that there exists $\mathcal{C}_{B}>0$ with 
$$\sup_{\epsilon\in D(0,\epsilon_0)}\left\|\tilde{b}_{jk}\right\|_{(\beta,\mu)}\le\mathcal{C}_B$$ for all $j,k\in\{0,1\}$ with $(j,k)\neq (0,1)$, and $\mathcal{C}_{C}>0$ such that
$$\sup_{\epsilon\in D(0,\epsilon_0)}\left\|C_{\ell}(m,\epsilon)\right\|_{(\beta,\mu)}\le\mathcal{C}_C,$$
for every $1\le \ell\le D-1$. We define
$$c_{\ell}(z,\epsilon)=\frac{1}{(2\pi)^{1/2}}\int_{-\infty}^{\infty}C_{\ell}(m,\epsilon)e^{izm}dm,\quad 1\le \ell\le D-1$$
$$b_{jk}(z,\epsilon)=\frac{1}{(2\pi)^{1/2}}\int_{-\infty}^{\infty}\tilde{b}_{jk}(m,\epsilon)e^{izm}dm,\quad j,k\in\{0,1\}, (j,k)\neq (0,1),$$
as in the general framework of the problem.

On the other hand, the geometric assumption on the choice of $d\in\R$ and the need of small ratio $d_D/k$ needed for Assumption (D) are no longer needed in this framework. Therefore, a wider choice in the geometry of the problem is available, together with a larger family of parameters involved in the problem.

The strategy to solve the problem is the same as that detailed in Section~\ref{secproblemstrategy}, by writing the solution of (\ref{eprals}) in the form
$$u(t,z,\epsilon)=u_0(t,z,\epsilon)+u_1(t,z,\epsilon)\frac{\log(\epsilon t)}{\log(q)},$$
for certain functions $u_0,u_1$ with $u_j(t,z,\epsilon):=U_j(\epsilon t, z,\epsilon)$, for $j\in\{0,1\}$ and $U_j$ given by the inverse Fourier transform and $q-$Laplace transform of adequate functions, i.e. 
$$u_{j}(t,z,\epsilon)=\frac{1}{(2\pi)^{1/2}}\frac{k}{\log(q)}\int_{-\infty}^{\infty}\int_{L_d}\frac{\omega_{j}(u,m,\epsilon)}{\Theta_{q^{1/k}}\left(\frac{u}{\epsilon t}\right)}\frac{du}{u}\exp(izm)dm,$$
along some direction $d$. From analogous steps as in Section~\ref{secproblemstrategy}, the problem is reduced to searching solutions of the following system of convolution equations:

\begin{multline}
Q(im)\omega_0(\tau,m,\epsilon)=\frac{\tau^{d_D}R_{D}(im)}{(q^{\frac{1}{k}})^{\frac{d_D(d_D-1)}{2}}}\omega_0(\tau,m,\epsilon)+\frac{d_D}{k}R_{D}(im)\frac{\tau^{d_D}}{(q^{\frac{1}{k}})^{\frac{d_D(d_D-1)}{2}}}\omega_1(\tau,m,\epsilon)\\
+\sum_{\ell=1}^{D-1}\epsilon^{\Delta_\ell-d_{\ell}}\frac{1}{(2\pi)^{1/2}}\int_{-\infty}^{+\infty}C_{\ell}(m-m_1,\epsilon)R_{\ell}(im_1)\left[\frac{\tau^{d_\ell}}{(q^{\frac{1}{k}})^{\frac{d_{\ell}(d_{\ell}-1)}{2}}}\sigma_{q;\tau}^{\delta_{\ell}-\frac{d_\ell}{k}}\omega_0( \tau,m_1,\epsilon)\right.\\
\hfill\left.+\delta_{\ell}\frac{\tau^{d_\ell}}{(q^{\frac{1}{k}})^{\frac{d_{\ell}(d_{\ell}-1)}{2}}}\sigma_{q;\tau}^{\delta_{\ell}-\frac{d_\ell}{k}}\omega_1( \tau,m_1,\epsilon)\right]dm_1\\
+\tilde{F}_0(\tau,m,\epsilon)+\frac{1}{(2\pi)^{1/2}}\int_{-\infty}^{+\infty}\tilde{b}_{00}(m-m_1,\epsilon)\omega_0(\tau,m_1,\epsilon)+\tilde{b}_{10}(m-m_1,\epsilon)\omega_1(\tau,m_1,\epsilon)dm_1\label{e376s}
\end{multline}

\begin{multline}
Q(im)\omega_1(\tau,m,\epsilon)=\frac{\tau^{d_D}R_{D}(im)}{(q^{\frac{1}{k}})^{\frac{d_D(d_D-1)}{2}}}\omega_1(\tau,m,\epsilon)\\
+\sum_{\ell=1}^{D-1}\epsilon^{\Delta_\ell-d_{\ell}}\frac{1}{(2\pi)^{1/2}}\int_{-\infty}^{+\infty}C_{\ell}(m-m_1,\epsilon)R_{\ell}(im_1)\frac{\tau^{d_\ell}}{(q^{\frac{1}{k}})^{\frac{d_{\ell}(d_{\ell}-1)}{2}}}\sigma_{q;\tau}^{\delta_{\ell}-\frac{d_\ell}{k}}\omega_1( \tau,m_1,\epsilon)dm_1\\
\hfill+\tilde{F}_1(\tau,m,\epsilon)+\frac{1}{(2\pi)^{1/2}}\int_{-\infty}^{\infty}\tilde{b}_{11}(m-m_1,\epsilon)\omega_1(\tau,m_1,\epsilon)dm_1\label{e377s}
\end{multline}

Observe that equation (\ref{e376s}) remains unchanged with respect to (\ref{e376}) whereas equation (\ref{e377s}) does not depend on $\omega_0$ when compared with (\ref{e377}). This is the key point of the simplified problem, as the system (\ref{e376s}), (\ref{e377s}) is triangular. One can proceed by solving (\ref{e377s}) first and obtain $\omega_1$, and then solve (\ref{e376s}) in $\omega_0$.

\subsection{Analytic solution of the auxiliary system}

The geometry of the problem needs less restrictive conditions, as we proceed to show. We choose 
$0<\rho<\min\{1,\frac{1}{2}q^{(d_D-1)/2k}\mathfrak{D}_1^{1/d_D}\}$ as in the general situation, and we choose $d\in\R$ such that $S_d\cup D(0,\rho)$ circumvents all the roots of $P_m(\tau)$ for every $m\in\R$. We observe that the statements of Lemma~\ref{lema477} hold. The number $\delta>0$ is chosen to satisfy $\hbox{dist}(S_d\cup \overline{D}(0,\rho),-\delta)\ge 1$.

In the following, we fix $\epsilon\in D(0,\epsilon_0)\setminus\{0\}$ and $\alpha > 0$ satisfying (\ref{e753}). We remark that Proposition~\ref{prop536} regarding the operator $\mathcal{H}_{\ell}$, and Proposition~\ref{prop598} concerning the forcing term are still valid in this situation. Proposition~\ref{prop565} is also valid, but in this framework a deep control of the norm of the operator $\mathcal{H}_P$ is not needed, which entangles more freedom on the coefficients and the geometry of the problem.

At this point, we define the operator
$$\mathcal{H}_{\epsilon}^{(1)}:\hbox{Exp}_{(k,\beta,\mu,\alpha,\rho)}^{q}\to \hbox{Exp}_{(k,\beta,\mu,\alpha,\rho)}^{q}$$
 by
$$\mathcal{H}_{\epsilon}^{(1)}(\omega_1(\tau,m))=\sum_{\ell=1}^{D-1}\epsilon^{\Delta_\ell-d_{\ell}}\mathcal{H}_{\ell}(\omega_1)+\frac{\tilde{F}_1(\tau,m,\epsilon)}{P_m(\tau)}+\frac{1}{P_m(\tau)}(\tilde{b}_{11}\star\omega_1).$$

\begin{prop}\label{prop641s}
There exists $\epsilon_0,\varpi_1,\varsigma_b>0$  such that for every $\epsilon\in D(0,\epsilon_0)\setminus\{0\}$ and $\mathcal{C}_{B}\le \varsigma_b$ (recall that $\mathcal{C}_{B}$ is determined by (\ref{e321})), the operator $\mathcal{H}_{\epsilon}^{(1)}$ satisfies that $\mathcal{H}_{\epsilon}^{(1)}(\overline{B}(0,\varpi_1))\subseteq \overline{B}(0,\varpi_1),$ where $\overline{B}(0,\varpi_1)\subseteq \hbox{Exp}_{(k,\beta,\mu,\alpha,\rho)}^{q}$ stands for the closed disc of radius $\varpi_1$ centered at the origin in $\hbox{Exp}_{(k,\beta,\mu,\alpha,\rho)}^{q}$. For every $\omega_{11},\omega_{12}\in \overline{B}(0,\varpi_1)$ it holds that
$$\left\|\mathcal{H}_{\epsilon}^{(1)}(\omega_{11}(\tau,m))-\mathcal{H}_{\epsilon}^{(1)}(\omega_{12}(\tau,m))\right\|_{(k,\beta,\mu,\alpha,\rho)}\le \frac{1}{2}\left\|\omega_{11}(\tau,m)-\omega_{12}(\tau,m)\right\|_{(k,\beta,\mu,\alpha,\rho)}.$$
\end{prop}
\begin{proof}
Let $\varpi_1>0$ such that $\tilde{C}_{F}\le \varpi_1/2$ and assume that $\epsilon_0,\varsigma_b>0$ are such that 
$$\sum_{\ell=1}^{D-1}\epsilon_0^{\Delta_{\ell}-d_\ell}C_{3,\ell}+\sup_{m\in\R}\left(\frac{1}{|Q(im)|}\right)\frac{1}{C_D}\varsigma_bC_2\le\frac{1}{2},$$
(recall $C_2$ is determined in Corollary~\ref{coro486}).
An analogous proof as that for (\ref{e648}) and (\ref{e687}) yields the result.
\end{proof}

As a direct consequence of the previous result, the operator $\mathcal{H}^{(1)}_{\epsilon}$ is contractive in $\overline{B}(0,\varpi_1)\subseteq \hbox{Exp}_{(k,\beta,\mu,\alpha,\rho)}^{q}$. Therefore, there exists a fixed point $\omega_1\in\hbox{Exp}_{(k,\beta,\mu,\alpha,\rho)}^{q}$ for such operator, with $\left\|\omega_1\right\|_{(k,\beta,\mu,\alpha,\rho)}\le\varpi_1$.

We observe from the construction that the solution depends holomorphically with respect to $\epsilon\in D(0,\epsilon_0)\setminus\{0\}$, as in the general situation. We define the function $(\tau,m,\epsilon)\mapsto \omega_{1}(\tau,m,\epsilon)$, where $\omega_{1}$ is the fixed point for $\mathcal{H}_{\epsilon}^{(1)}$, and adopt the same notation for such function. We recall that
\begin{equation}\label{e1231}
\sup_{\epsilon\in D(0,\epsilon_0)\setminus\{0\}}|\omega_{1}(\tau,m,\epsilon)|\le \varpi_1 \frac{1}{(1+|m|)^{\mu}}e^{-\beta |m|}\exp\left(\frac{k}{2}\frac{\log^2|\tau+\delta|}{\log(q)}+\alpha\log|\tau+\delta|\right),
\end{equation}
for every $\tau\in S_{d}\cup\overline{D}(0,\rho)$, $m\in\R$ and $\epsilon\in D(0,\epsilon_0)\setminus\{0\}$.

Let $\epsilon\in D(0,\epsilon_0)\setminus\{0\}$. We define the operator
$$\mathcal{H}_{\epsilon}^{(0)}:\hbox{Exp}_{(k,\beta,\mu,\alpha,\rho)}^{q}\to \hbox{Exp}_{(k,\beta,\mu,\alpha,\rho)}^{q}$$
 by
$$\mathcal{H}_{\epsilon}^{(0)}(\omega_0(\tau,m))=\sum_{\ell=1}^{D-1}\epsilon^{\Delta_\ell-d_{\ell}}\mathcal{H}_{\ell}(\omega_0)+\frac{1}{P_m(\tau)}(\tilde{b}_{00}\star\omega_0)+G_{\epsilon}(\tau,m),$$
where
$$G_{\epsilon}(\tau,m):=\mathcal{H}_P(\omega_1)+\sum_{\ell=1}^{D-1}\epsilon^{\Delta_\ell-d_{\ell}}\delta_{\ell}\mathcal{H}_{\ell}(\omega_1)+\frac{\tilde{F}_0(\tau,m,\epsilon)}{P_m(\tau)}+\frac{1}{P_m(\tau)}(\tilde{b}_{10}\star\omega_1).$$

\begin{lemma}
The function $G_{\epsilon}$ belongs to $\hbox{Exp}_{(k,\beta,\mu,\alpha,\rho)}^{q}$ for every $\epsilon\in D(0,\epsilon_0) \setminus \{ 0 \}$. Moreover, for every $\tau\in S_d\cup \overline{D}(0,\rho)$ and $m\in\R$, the map $\epsilon\mapsto G_{\epsilon}(\tau,m)$ is holomorphic on $D(0,\epsilon_0)\setminus\{0\}$.
\end{lemma}
\begin{proof}
Recall that for every $\epsilon\in D(0,\epsilon_0)\setminus\{0\}$, the function $(\tau,m) \mapsto \omega_1(\tau,m,\epsilon)$ belongs to $\hbox{Exp}_{(k,\beta,\mu,\alpha,\rho)}^{q}$. Therefore, a direct application of Proposition~\ref{prop536}, Proposition~\ref{prop565}, Proposition~\ref{prop598}, and Corollary~\ref{coro486} together with Lemma~\ref{lema425} allows us to conclude that $G_{\epsilon}(\tau,m)\in \hbox{Exp}_{(k,\beta,\mu,\alpha,\rho)}^{q}$. The holomorphic dependence with respect to the perturbation parameter is a consequence of the construction involved in such results.
\end{proof}

\begin{prop}\label{prop641ss}
There exists $\epsilon_0,\varpi_0,\varsigma_b>0$  such that for every $\epsilon\in D(0,\epsilon_0)\setminus\{0\}$ and $\mathcal{C}_{B}\le \varsigma_b$, the operator $\mathcal{H}_{\epsilon}^{(0)}$ satisfies that $\mathcal{H}_{\epsilon}^{(0)}(\overline{B}(0,\varpi_0))\subseteq \overline{B}(0,\varpi_0),$ where $\overline{B}(0,\varpi_0)\subseteq \hbox{Exp}_{(k,\beta,\mu,\alpha,\rho)}^{q}$ stands for the closed disc of radius $\varpi_0$ centered at the origin in $\hbox{Exp}_{(k,\beta,\mu,\alpha,\rho)}^{q}$. For every $\omega_{01},\omega_{02}\in \overline{B}(0,\varpi_0)$ it holds that
$$\left\|\mathcal{H}_{\epsilon}^{(0)}(\omega_{01}(\tau,m))-\mathcal{H}_{\epsilon}^{(0)}(\omega_{02}(\tau,m))\right\|_{(k,\beta,\mu,\alpha,\rho)}\le \frac{1}{2}\left\|\omega_{01}(\tau,m)-\omega_{02}(\tau,m)\right\|_{(k,\beta,\mu,\alpha,\rho)}.$$
\end{prop}
\begin{proof}
Let $\varpi_0>0$ large enough in order that 
$$\left[\frac{d_D}{k}\frac{1}{\mathfrak{D}_1}\max\{\frac{1}{C_D},\frac{1}{\mathfrak{D}_3}\}+\sum_{\ell=1}^{D-1}\epsilon_0^{\Delta_\ell-d_\ell}\delta_\ell C_{3,\ell}+\frac{1}{C_D}\max_{m\in\R}\left(\frac{1}{|Q(im)|}\right)\varsigma_bC_2\right]\varpi_1+\tilde{C}_{F}\le\frac{\varpi_0}{2},$$
where $\varpi_1$ is determined in Proposition~\ref{prop641s}, $C_2$ is fixed in Corollary~\ref{coro486}, and $\varsigma_b,\epsilon_0,C_2>0$ are such that
$$\sum_{\ell=1}^{D-1}\epsilon_0^{\Delta_\ell-d_{\ell}} C_{3,\ell}+\sup_{m\in\R}\left(\frac{1}{|Q(im)|}\right)\frac{1}{C_D}\varsigma_bC_2\le\frac{1}{2}.$$
Then, in view of Proposition~\ref{prop536}, Proposition~\ref{prop565}, Corollary~\ref{coro486} and Proposition~\ref{prop598}, we have that
\begin{multline*}
\left\|G_{\epsilon}\right\|_{(k,\beta,\mu,\alpha,\rho)}\le \left\|\mathcal{H}_P(\omega_1)\right\|_{(k,\beta,\mu,\alpha,\rho)}+\sum_{\ell=1}^{D-1}\epsilon_0^{\Delta_\ell-d_{\ell}}\delta_{\ell}\left\|\mathcal{H}_{\ell}(\omega_1)\right\|_{(k,\beta,\mu,\alpha,\rho)}\\
\hfill+\left\|\frac{\tilde{F}_0(\tau,m,\epsilon)}{P_m(\tau)}\right\|_{(k,\beta,\mu,\alpha,\rho)}+\left\|\frac{1}{P_m(\tau)}(\tilde{b}_{10}\star\omega_1)\right\|_{(k,\beta,\mu,\alpha,\rho)}\\
\le \left(\frac{d_D}{k}\frac{1}{\mathfrak{D}_1}\max\{\frac{1}{C_D},\frac{1}{\mathfrak{D}_3}\}+\sum_{\ell=1}^{D-1}\epsilon_0^{\Delta_\ell-d_\ell}\delta_\ell C_{3,\ell}+\frac{1}{C_D}\max_{m\in\R}\left(\frac{1}{|Q(im)|}\right)\mathcal{C}_BC_2\right)\varpi_1+ \tilde{C}_F\\
\le \left(\frac{d_D}{k}\frac{1}{\mathfrak{D}_1}\max\{\frac{1}{C_D},\frac{1}{\mathfrak{D}_3}\}+\sum_{\ell=1}^{D-1}\epsilon_0^{\Delta_\ell-d_\ell}\delta_\ell C_{3,\ell}+\frac{1}{C_D}\max_{m\in\R}\left(\frac{1}{|Q(im)|}\right)\varsigma_bC_2\right)\varpi_1+ \tilde{C}_F\le \varpi_0/2,
\end{multline*}
for every $\epsilon\in D(0,\epsilon_0)\setminus\{0\}$.
Take $\omega_0\in \overline{B}(0,\varpi_0)$. An analogous argument as before allows us to arrive at
\begin{multline*}
\left\|\mathcal{H}_{\epsilon}^{(0)}(\omega_0(\tau,m))\right\|_{(k,\beta,\mu,\alpha,\rho)}\hfill\\
\le\sum_{\ell=1}^{D-1}\epsilon_0^{\Delta_\ell-d_{\ell}}\left\|\mathcal{H}_{\ell}(\omega_0)\right\|_{(k,\beta,\mu,\alpha,\rho)}+\left\|\frac{1}{P_m(\tau)}(\tilde{b}_{00}\star\omega_0)\right\|_{(k,\beta,\mu,\alpha,\rho)}+\left\|G_{\epsilon}(\tau,m)\right\|_{(k,\beta,\mu,\alpha,\rho)}\\
\le \left(\sum_{\ell=1}^{D-1}\epsilon_0^{\Delta_\ell-d_{\ell}} C_{3,\ell}+\sup_{m\in\R}\left(\frac{1}{|Q(im)|}\right)\frac{1}{C_D}\mathcal{C}_BC_2\right)\varpi_0+\frac{\varpi_0}{2}\\
\le \left(\sum_{\ell=1}^{D-1}\epsilon_0^{\Delta_\ell-d_{\ell}} C_{3,\ell}+\sup_{m\in\R}\left(\frac{1}{|Q(im)|}\right)\frac{1}{C_D}\varsigma_bC_2\right)\varpi_0+\frac{\varpi_0}{2}\le \frac{\varpi_0}{2}+\frac{\varpi_0}{2}=\varpi_0.
\end{multline*}
The second part of the proof follows analogous arguments as in the previous results.
\end{proof}

The operator $\mathcal{H}^{(0)}_{\epsilon}$ is contractive in $\overline{B}(0,\varpi_0)\subseteq \hbox{Exp}_{(k,\beta,\mu,\alpha,\rho)}^{q}$. Therefore, there exists a fixed point $\omega_0\in\hbox{Exp}_{(k,\beta,\mu,\alpha,\rho)}^{q}$ for such operator, with $\left\|\omega_0\right\|_{(k,\beta,\mu,\alpha,\rho)}\le\varpi_0$. The assignment of $\omega_{0}$ departing from $\epsilon$ can be done holomorphically with respect to $\epsilon\in D(0,\epsilon_0)\setminus\{0\}$, following the same reasoning as in the general case. We define the function $(\tau,m,\epsilon)\mapsto \omega_{0}(\tau,m,\epsilon)$, where $\omega_{0}$ is the fixed point for $\mathcal{H}_{\epsilon}^{(0)}$, and adopt the notation $\omega_0$ for such function, as above. We have
\begin{equation}\label{e1231b}
\sup_{\epsilon\in D(0,\epsilon_0)\setminus\{0\}}|\omega_{0}(\tau,m,\epsilon)|\le \varpi_0 \frac{1}{(1+|m|)^{\mu}}e^{-\beta |m|}\exp\left(\frac{k}{2}\frac{\log^2|\tau+\delta|}{\log(q)}+\alpha\log|\tau+\delta|\right),
\end{equation}
for every $\tau\in S_{d}\cup\overline{D}(0,\rho)$, $m\in\R$ and $\epsilon\in D(0,\epsilon_0)\setminus\{0\}$.

At this point, one can derive the existence of a solution to equation (\ref{eprals}). Indeed, let us fix a good covering $(\mathcal{E}_{p})_{0\le p\le \zeta-1}$ in $\C^{\star}$ and let $\{(\mathcal{R}_{d_p,\Delta})_{0\le p\le \zeta-1}, D(0,\rho),\mathcal{T}\}$ be a set associated to the previous good covering, where $\zeta\ge 2$ and with $\Delta,\rho$ and $d_p$ for $0\le p\le \zeta-1$ as in Section~\ref{sec9}. We choose $\alpha$ as in (\ref{e753}).

\begin{theo}\label{teo5}
In the previous situation, for every $0\le p\le \zeta-1$, the equation (\ref{eprals}) admits a solution $u_p(t,z,\epsilon)$ of the form
\begin{equation}\label{e1272}
u_p(t,z,\epsilon)=u_{0,p}(t,z,\epsilon)+u_{1,p}(t,z,\epsilon)\frac{\log(\epsilon t)}{\log(q)},
\end{equation}
holomorphic on $\mathcal{T}\times H_{\beta'}\times \mathcal{E}_p$, for every $0<\beta'<\beta$. Moreover, there exist constants $\tilde{K}>0$ and $\tilde{\alpha}\in\R$ such that
\begin{equation}\label{e1276}
\sup_{t\in\mathcal{T},z\in H_{\beta'}}|u_{j,p+1}(t,z,\epsilon)-u_{j,p}(t,z,\epsilon)|\le \tilde{K}\exp\left(-\frac{k}{2\log(q)}\log^2|\epsilon|\right)|\epsilon|^{\tilde{\alpha}} 
\end{equation}
for $\epsilon\in \mathcal{E}_{p}\cap\mathcal{E}_{p+1}$, and all $0\le p\le \zeta-1$ (by identifying $u_{j,\zeta}$ with $u_{j,0}$), for $j\in\{0,1\}$.
\end{theo}
\begin{proof}
Let $\varpi_1$ be determined in Proposition~\ref{prop641s} and $\varpi_0$ as in Proposition~\ref{prop641ss}. We choose $\epsilon_0>0$ and $\varsigma_b>0$ small enough which guarantee the existence of a fixed point for $\mathcal{H}_\epsilon^{(1)}$, say $\omega_{1,p}\in\overline{B}(0,\varpi_1)\subseteq\hbox{Exp}_{(k,\beta,\mu,\alpha,\rho)}^{q}$. In view of Proposition~\ref{prop641s}, $\mathcal{H}_{\epsilon}^{(1)}$ is a contractive map from a closed disc contained in a Banach space into itself, which depends holomorphically on $\epsilon$. The estimates (\ref{e1231}) allow us to construct 
$$u_{1,p}(t,z,\epsilon)=\frac{1}{(2\pi)^{1/2}}\frac{k}{\log(q)}\int_{-\infty}^{\infty}\int_{L_d}\frac{\omega_{1,p}(u,m,\epsilon)}{\Theta_{q^{1/k}}\left(\frac{u}{\epsilon t}\right)}\frac{du}{u}\exp(izm)dm,$$
which is holomorphic on $\mathcal{T}\times H_{\beta'}\times \mathcal{E}_p$, for any fixed $0<\beta'<\beta$ as in Section~\ref{sec9}. Finally, a second fixed point argument applied to $\mathcal{H}^{(0)}_{\epsilon}$ guarantees the existence of a fixed point for $\mathcal{H}_\epsilon^{(0)}$, say $\omega_{0,p}\in\overline{B}(0,\varpi_0)\subseteq\hbox{Exp}_{(k,\beta,\mu,\alpha,\rho)}^{q}$. After adopting the same notation with respect to $\epsilon$, one can define the function
$$u_{0,p}(t,z,\epsilon)=\frac{1}{(2\pi)^{1/2}}\frac{k}{\log(q)}\int_{-\infty}^{\infty}\int_{L_d}\frac{\omega_{0,p}(u,m,\epsilon)}{\Theta_{q^{1/k}}\left(\frac{u}{\epsilon t}\right)}\frac{du}{u}\exp(izm)dm,$$
which is holomorphic on $\mathcal{T}\times H_{\beta'}\times \mathcal{E}_p$, for any fixed $0<\beta'<\beta$. The function (\ref{e1272}) turns out to be a solution of (\ref{eprals}) from the properties stated in Section~\ref{secreview} on the integral transforms.

The statement (\ref{e1276}) follows from a word-by-word proof with respect to that of Theorem~\ref{teodif}.
\end{proof}

The same steps allow us to conclude the main result in this particularized framework. The proof is omitted. As in the main result of the present work, Theorem~\ref{teopral}, we fix $0<\beta'<\beta$ and consider be the Banach space of bounded holomorphic functions on $\mathcal{T}\times H_{\beta'}$, with the norm of the supremum, denoted by $\mathbb{E}$.

\begin{theo}
There exist $\hat{u}_0,\hat{u}_1\in\mathbb{E}[[\epsilon]]$ such that the function $\epsilon\mapsto u_{j,p}$, constructed in Theorem~\ref{teo5}, admits $\hat{u}_j$ as its $q-$Gevrey asymptotic expansion of order $1/k$ on $\mathcal{E}_p$ for every $0\le p\le \zeta-1$ and $j\in\{0,1\}$.

In addition to this, the formal expression
\begin{equation}\label{e974b}
\hat{u}(t,z,\epsilon)=\hat{u}_0(t,z,\epsilon)+\hat{u}_1(t,z,\epsilon)\frac{\log(\epsilon t)}{\log(q)}
\end{equation}
is a formal solution of
\begin{multline*}
Q(\partial_z)\hat{u}(t,z,\epsilon)=(\epsilon t)^{d_D}\sigma_{q;t}^{\frac{d_D}{k}}R_D(\partial_{z})\hat{u}(t,z,\epsilon)+\sum_{\ell=1}^{D-1}\epsilon^{\Delta_\ell}t^{d_{\ell}}c_{\ell}(z,\epsilon)R_{\ell}(\partial_z)\sigma_{q;t}^{\delta_{\ell}}\hat{u}(t,z,\epsilon)\\
+f(t,z,\epsilon)+b_{00}(z,\epsilon)\left[\hat{u}(t,z,\epsilon)-\left(\frac{1}{2\pi i}(\gamma_{\epsilon}^{*}-\hbox{id})\hat{u}(t,z,\epsilon)\right)\log(\epsilon t)\right]\\
+\left(b_{10}(z,\epsilon)+b_{11}(z,\epsilon)\frac{\log(\epsilon t)}{\log(q)}\right)\left[\frac{\log(q)}{2\pi i}(\gamma_{\epsilon}^{*}-\hbox{id})\hat{u}(t,z,\epsilon)\right].
\end{multline*}
with $\gamma_{\epsilon}^{*}$ being the formal monodromy operator defined in (\ref{e125}).
\end{theo}

\end{document}